\mag=\magstep1
\hsize 4.8in
\vsize 7.6in
\input epsf.sty

\def\qed{\hfill $ \sqcup\!\!\!\!\sqcap$\medskip}
\def\ni{\noindent}

\rightline{November 2000}\medskip

\centerline{\bf Sextactic points on a simple closed curve}

\medskip
\vskip .7cm
\centerline{\it by Gudlaugur Thorbergsson and Masaaki Umehara}
\vskip 1.1cm
\bigskip

{ \ni{\bf Abstract.} We give optimal lower bounds for the number of
sextactic
points on a simple closed curve in the real projective plane. Sextactic
points are after inflection points the simplest projectively invariant
singularities on such curves. Our method is axiomatic and can be applied in
other situations.
\par}\bigskip

{\ni{\bf Content}\medskip}
{ \ni 1. Introduction\par}
{ \ni 2. Preliminaries\par}
{ \ni 3. Intrinsic conic systems\par}
{ \ni 4. An application to strictly convex curves\par}
{ \ni 5. An application to simple closed curves\par}
{ \ni Appendix A: Simple closed curves with few inflection points\par}
{ \ni Appendix B: Examples of curves with few sextactic points\par}
{ \ni  Appendix C: Sextactic points on a complex plane algebraic curve\par}

\bigskip
\noindent{\bf \S 1\hskip 0.1in Introduction}\medskip

In analogy with tangent lines and inflection points of regular curves in
the real or complex projective plane, one can consider their osculating
conics and
sextactic points. Choose five points on a curve $\gamma$ in a neighborhood
of a point $p$ on $\gamma$
that is not an inflection point. There is a unique regular conic passing
through the five points. Letting the five points all converge to $p$, the
conics converge
to a uniquely defined regular conic that is called the {\it osculating
conic of $\gamma$ in $p$}. The
osculating conic meets $\gamma$ with multiplicity at least five in $p$. If
it meets with multiplicity
at least six in $p$, then $p$ is called a {\it sextactic point.} \medskip

Inflection and sextactic points on curves in the complex projective plane
were well
understood already in the nineteenth century. We will make some historic
remarks
on this towards the end of
the introduction. It is the case of curves in the real projective plane
that still poses
problems.\medskip

In the present paper we will be dealing with closed $C^\infty$-parameterized
curves
$\gamma:S^1\to P^2$ that are simple (free of self-intersections) and
regular (nowhere vanishing tangent
vector). Here and elsewhere in the paper we let $P^2$ denote the {\it real}
projective plane.
The existence of inflection and sextactic points on such curves has of
course been much
studied. Of importance for us  is the result of M\"obius [M\"o] that a
simple regular curve in
$P^2$ that is not null-homotopic has at least three inflection points.
As far as we know,  the first
paper to deal with sextactic points on curves in the real projective (or
affine) plane that are not necessarily algebraic,  is the
paper [Mu 1] of Mukhopadhyaya from the year 1909. There it is proved that a
strictly convex curve in
the affine plane has at least six sextactic points. An independent proof of
this theorem due to
Herglotz and Radon  was published by Blaschke [Bl 1] in 1917.
Proofs can also be found in the textbooks [Bl 2] and [Bo].
In [Mu 2] Mukhopadhyaya proved that three of
these sextactic points can be chosen so that the corresponding osculating
conics are inscribed and
another three such that the corresponding osculating conics are
circumscribed. We are not aware of
any results on sextactic points on curves that are not strictly
convex. For  recent papers on sextactic points
on strictly convex curves and related matters, see e.g. [Ar
3] and [GMO].
\medskip

Our main results are summarized  in the following three theorems. We have
not tried to
state here everything in its strongest form. More precise results can be
found in sections four and
five. Notice that we will give examples in Appendices B and C showing that
all these theorems are
optimal. Two of the examples in Appendix B were communicated to us by
Izumiya and Sano [IS] who
found them  in their study of affine evolutes.\medskip

\noindent{\bf 1.1. Theorem.} {\it Let $\gamma$ be a simple closed curve in
$P^2$ that is not nullhomotopic. Then $\gamma$ has at least three sextactic
points.}\medskip

The result of M\"obius mentioned above is one of the essential ingredients
in the
proof of this theorem. Notice that the theorem is optimal since the
noncontractible branch of a real
cubic has exactly three sextactic points as we will explain in Appendix C.
Notice also that the theorem was stated as a problem by Bol in [Bo] on
p.~43. A sketch of 
a proof of Theorem 1.1 under rather strong genericity assumptions on the
inflection points of $\gamma$ was given by Fabricius-Bjerre in [Fa]; see
Remark (iii) after Proposition 5.1.\medskip

\noindent{\bf 1.2. Theorem.} {\it Let $\gamma$ be a
simple closed curve in $P^2$ that is nullhomotopic.
\item{(i)} (Mukhopadhyaya) If $\gamma$ is strictly convex, then it has at
least six
sextactic points.
\item{(ii)} If $\gamma$ is not convex, then it has at least three sextactic
points.
\item{(iii)} If $\gamma$ is convex, then it has at least two sextactic
points.}\medskip

Part (i), or Mukhopadhyaya's theorem, is optimal, since a nullhomotopic
component of a regular real
cubic is strictly convex and has exactly six sextactic points as we will
explain in Appendix C.
That the other parts are optimal will be explained in Appendix B.
\medskip

Counting sextactic points and inflection points together, we can prove the
following
theorem.\medskip

\noindent{\bf 1.3. Theorem.} {\it Let $\gamma$ be a
simple closed curve in $P^2$ that is nullhomotopic. Then the total number of
sextactic and inflection points on $\gamma$ is at least four}.\medskip

This theorem is optimal as an example of Izumiya and Sano [IS] shows that
we explain in Appendix
B. It will be clear from the
proof that Theorem 1.3 can
only be optimal for convex curves
with one or two inflection points.\medskip

The arguments in  our proofs are inspired by  those of
Mukhopadhyaya in [Mu 1] and especially in [Mu 2], although there are of
course
new ideas needed to deal with curves with inflection points. We have chosen
an axiomatic approach that is similar in spirit to the one introduced by
the second author in [Um] to deal with vertices and was further studied in
[TU 1].
The main idea behind this approach was motivated by the paper [Kn] of H.
Kneser.
It should be pointed out that our theorems are more generally true for
curves that are only
$C^4$  with essentially the same proofs, see the remark after Proposition
5.1, but we stay in
the $C^\infty$-category to simplify the exposition. Notice that one has to
modify the
definition of a sextactic point in the case of $C^4$-curves, see section
two.
\medskip

We would like to make a few remarks on inflection and sextactic points on
algebraic curves in the
complex projective plane.
There is a formula due to Pl\"ucker (1835) that one can find in most
textbooks on algebraic curves
saying that a regular algebraic curve
$\gamma$ of degree
$d$ in
$P^2({\bf C})$ has exactly
$3d(d-2)$ inflection points counted with multiplicities. It is much less
known that Cayley [Ca 2]
proved  in 1865 that such a curve (with simple inflection points) has
exactly $3d(4d-9)$ sextactic
points counted with multiplicities.
The condition which we have put within parentheses is not in Cayley's paper
although it is needed as
we will see in Appendix C.
Pl\"ucker and Cayley used the same strategy of proof: there is a curve of
degree $3(d-2)$
that intersects $\gamma$ precisely in the inflection points and similarly
there is a curve of degree
$3(4d-9)$ that intersects $\gamma$ precisely in the sextactic points. The
results then follow from
B\'ezout's theorem. \medskip

 The term {\it sextactic point} might have been introduced by Cayley in [Ca
1].
Cayley remarks that sextactic points were studied before him by Pl\"ucker
and Steiner without
giving concrete references. He is certainly referring to papers in Crelle's
Journal {\bf 32} (1846)
by Steiner and {\bf 34} (1847) by Pl\"ucker. One can add a paper by Hesse
in volume {\bf 36} (1848)
of the same journal. In all of these papers it is claimed that there are
twenty seven sextactic
points on a (smooth) cubic. Steiner claims in his paper that does not
contain any proofs that nine
of these are always real. This is only correct as Pl\"ucker points out if
the curve has two real
branches. A real cubic has three sextactic points if it has only one real
branch. We will discuss
this in Appendix C. Pl\"ucker's paper is a
polemic against Steiner and his methods in favor of analytic
geometry.\medskip

A formula due to Klein implies that a smooth algebraic curve of degree $d$
in the real projective
plane can have at most $d(d-2)$ inflection points, i.e., at most one third
of the complex
inflection points can be real. An
analogous result for
sextactic points seems to be unknown. A rigorous proof of Klein's formula
was given by Wall in
[Wa].
\medskip

The content of the sections of the paper is as follows.
Section two contains preliminaries. Section three explains our axiomatic
approach
to sextactic points. In section four we give a complete proof of the
results of Mukhopadhyaya since we need all the arguments involved, and a
treatment
of these ideas satisfying modern standards does not seem to exist.
In section five we prove the above theorems (Theorem 1.1 is
the same as 5.2, Theorem 1.2 (ii) is in 5.3 and 5.5, (iii) is in 5.4,
Theorem 1.3 is in 5.4 and 5.5.). In
Appendix A we prove a theorem on simple closed curves in $P^2$ with few
inflection points that is
needed in section five. In Appendix B we give examples that show that the
above theorems and some
of the results in section five are optimal. Two of these examples are due to
Izumiya and Sano [IS].
In Appendix C we sketch a proof of
the theorem of Cayley mentioned above that is based on standard results on
inflection points of
linear systems. What we prove is slightly more general than Cayley's result
since we do not make any
assumptions on the multiplicity of the inflection points of the curve. We
also discuss the
sextactic points on cubics in Appendix C.

\bigskip
\noindent{\bf \S 2\hskip 0.1in Preliminaries}\medskip

\ni{\bf A. Multiplicity of intersection points}\medskip

Let
$\gamma$ and $\sigma$ be two smooth and regular parameterized
curves in $P^2$.  The following definitions are all of a local nature. We
therefore assume that both curves are {\it simple}, i.e., without
self-intersections. Assume that
$p\in P^2$ lies in the image of both curves.
 Then the  {\it multiplicity of the intersection of
$\gamma$ and
$\sigma$ in $p$} is defined as follows.  The multiplicity
is equal to one if
$\gamma$ and $\sigma$  intersect transversally in
$p$.  If they do not
intersect transversally, we look at coordinates
$(x,y)$ about $p$ that we
assume to correspond to $(0,0)$ with the
$x$-axis as the common
tangent. Express the curves locally as graphs
over the $x$-axis in these
coordinates. Assume that the first $k$
derivatives of the $y$-components
of the curves coincide in
$0$, but
not the $k+1$st. Then we say that the {\it multiplicity of the
common
point $p$
 is equal to $k+1$} and the {\it order of contact of
$\gamma$ and
$\sigma$ in
$p$ is equal to $k$}.  If all derivatives
coincide, the {\it multiplicity and
the order of contact are
infinite.}

We say that $\gamma$ and $\sigma$ {\it cross in $p$} if
 either they meet
transversally in $p$ or if $p$ is isolated in $\gamma\cap\sigma$ and
there are
coordinates $(x,y)$ in which $p$
corresponds to $(0,0)$, the tangent lines  of both curves
in $p$
corresponds to the $x$-axis,  the images of $\gamma$ and
$\sigma$ are locally around $p$
graphs of functions $f$ and $g$, and $f-g$ changes sign
in $0$ and only vanishes in $0$. We say
that they {\it are  locally one on the side of the
other around $p$}, if they are
not transversal in $p$, $f-g$ does not
change sign in $p$ and $f-g$ does
not vanish except in $0$ for
functions $f$ and $g$ as above.

If two curves $\gamma$ and $\sigma$
have an isolated connected
set
$J$ of
common points, we can extend the above definition and say that
the
curves either {\it cross in $J$} or {\it are locally one on the
side of the
other around $J$.}

If $\gamma$ and $\sigma$ meet
with finite multiplicity in a point
$p$, then it follows from
Elementary Calculus that $p$ is isolated in the
set of common points.
If
$\gamma$ and $\sigma$  meet with odd multiplicity in $p$, then
it
follows  that the curves  cross in
$p$.  If the multiplicity is
even, they are locally one on the side of each
other.\bigskip

\ni{\bf B. Conics}\medskip
 For us
a {\it conic} will be a quadric  without singularities. This
excludes
the reducible quadrics which are either a union of two lines
or a  line
counted twice.  Any two different conics are projectively
equivalent.
\medskip

 Two conics in $P^2$ are identical if they have
five different points in
common.  Given five points in $P^2$, no three
of which are collinear,
there is a unique conic passing through these
points.\medskip

 If we count  common points with multiplicities in
the sense defined
above, it follows that two conics with five points
in common coincide. The family of conics that are tangent to a
curve $\gamma$ at two
different points $p$ and $q$ is one dimensional,
and two conics in that
family only have the two given points in
common.
There is  a one dimensional family of conics
that  meet a curve $\gamma$ with multiplicity at least four at a given
point,
 and
two conics in that family only have the given
point in common. \medskip

\bigskip

\ni{\bf C. Inflection and sextactic points}\medskip

We will call a point $p$ on $\gamma$ an {\it  inflection point} if
$det(\hat\gamma,\hat\gamma^\prime,\hat\gamma^{\prime\prime})$
vanishes in $p$,
where $\hat\gamma$ is a representation of $\gamma$ in homogeneous
coordinates.  An
equivalent definition is to say that $p$ on $\gamma$ is an
inflection point if  $\gamma$ and the tangent line
of $\gamma$ in
$p$ meet with multiplicity at least three in $p$.
We are used to think of an inflection point as a point where, roughly
speaking,
 the direction changes in which the curve is bending.
We therefore call $p$ on $\gamma$
a {\it  true inflection point}, if
the tangent line of
$\gamma$ at $p$ and
$\gamma$ cross in $p$ or if they cross in the connected
component containing $p$ of their common points.
 It
follows that the multiplicity with which
 the tangent line of
$p$ at $\gamma$
and
$\gamma$ meet is odd or  infinite in $p$, if $p$ is a true inflection point.
\medskip

We can now state a more precise version of the Theorem of M\"obius [M\"o]
than in the
introduction, see [TU 1].\medskip

\ni{\bf 2.1. Theorem.} (M\"obius) {\it Let $\gamma$ be a simple closed
curve in $P^2$
that is not null-homotopic and not a projective line. Then $\gamma$ has at
least three intervals of
true inflection points.}\medskip

We will also need the following theorem. The complements of lines in $P^2$
are called {\it affine
planes}.\medskip

\ni{\bf 2.2. Theorem.} {\it Let $\gamma: [0,1]\to P^2$ be a regular simple
arc such that
no $\gamma(t)$ for $t\in (0,1)$ is an inflection point. Then $\gamma$ lies
in an affine plane.

A contractible simple closed curve $\gamma: S^1\to P^2$ with less than four
intervals of true
inflection points is contained in an affine plane.}\medskip

The first part of this theorem is in [Ar 2], see also [Ar 1].  The second
part of the theorem is
proved in [Ar 2]  for the special case of precisely two true inflection
points. We give a complete
proof of the second part of this proposition in Appendix A, that is based
on our paper [TU 1].\medskip

 A closed curve $\gamma$ in $P^2$ is called {\it convex} if
it lies in some
affine plane in $P^2$ where it is convex in the usual sense of bounding a
convex domain. A closed curve is called
{\it strictly convex} if it is convex and has no inflection points. One can
prove that a closed curve without self-intersections and inflection points
in $P^2$ is
a strictly convex
closed curve. In fact, Theorem 2.2
implies that
such a curve is contained in an affine plane,  where the claim is standard.
\medskip

 Let $p$ be a point on a smooth and regular
curve $\gamma$ in
$P^2$ that is not an inflection point.  Then there is
 a unique conic that meets $\gamma$ with multiplicity five at least in $p$,
see
e.g.~[Bo]. This conic is called the {\it
osculating conic} of $\gamma$ at
$p$. It is clear that
there is no regular conic meeting a curve with multiplicity five or higher
in an inflection point. If the multiplicity is precisely
five between a curve $\gamma$ and the osculating conic in $p$, then
$\gamma$ and the osculating conic cross in $p$.
\medskip

If the osculating conic of $\gamma$ at $p$ meets $\gamma$ with multiplicity
six at least in $p$, then $p$ is called a {\it sextactic point.} If
$\gamma$ lies in
an affine plane $A^2\subset P^2$, then $p$ is sextactic if and only if $p$
is a critical point of the affine curvature of $\gamma$, see [Bl 2]. The
affine
curvature (or the projective length element) will not play any role in the
proofs of the main
results of this paper and will only be referred to in some remarks.\medskip

We will need the following lemma that can already be found in [Mu 1]. The
books [Bl 2] and [Bo] bring it as an exercise. We will give
 a proof of the lemma in (4.9), which applies to
arcs that are only $C^4$.
\medskip

\ni{\bf 2.3. Lemma.} {\it Let $\gamma$ be an arc in $P^2$ that is free of
inflection and sextactic points. Then the osculating conics at two different
points of
$\gamma$ do not meet.}\medskip

Our methods will mostly imply the existence of sextactic points with the
property
that the curve and the osculating conic do not cross there. In fact,
Mukhopadhyaya
requires this noncrossing property in his definition of a sextactic point.
\medskip

We now
introduce  terminology to describe the different cases of sextactic points
we will
encounter. Notice that a conic divides $P^2$ into two
closed domains, one of which is a homeomorphic to a disk, the other is
homeomorphic to a
M\"obius strip. We say that a curve is {\it inside} the conic if it lies in
the disk and
{\it outside} if it lies in the M\"obius strip.
We will call a sextactic point $p$ of $\gamma$ {\it minimal} if
some arc of $\gamma$
around $p$ is
 inside of the osculating conic at $p$ and {\it maximal} if
 some arc around
$p$  is outside the osculating conic at $p$. We will call a sextactic point
of a closed curve
$\gamma$ {\it globally maximal} if the whole curve $\gamma$ lies inside the
osculating conic and
{\it globally minimal} if it lies outside the osculating conic. A sextactic
point of a closed
curve
$\gamma$ will be called {\it clean} if the intersection of the osculating
conic and $\gamma$
is connected. \medskip

Notice that a sextactic point in which $\gamma$ and the osculating conic
meet with odd
multiplicity does not satisfy these additional properties we have been
defining and the same
{\it can} happen if the multiplicity is infinite.\medskip

The above definition of a sextactic point only makes sense for curves that
are
$C^5$ at least. A point on a $C^4$-regular arc is called {\it sextactic} if
the
osculating conic does not cross in that point. Notice that this definition
implies, but is not equivalent to the original definition if the curve is
$C^5$.

\bigskip
\noindent
{\bf \S 3\hskip 0.1in Intrinsic conic systems} \medskip

In this section we explain our axiomatic approach to sextactic points.
It will be the main tool to prove the existence of sextactic points in the
 several different, although similar, situations in sections four and five.
We will define
an abstract notion of a sextactic point in our axiomatic setting that will
turn out to correspond to
those sextactic points of curves that we call maximal or minimal, see the
Preliminaries. \medskip

\ni{\bf A. Intrinsic circle systems}\medskip

We will need a lemma on intrinsic circle systems.
Let $I$ be either the circle $S^1$ or an interval of $S^1$ that can be
open, closed or
halfopen. We denote the closure of $I$ by $\bar I$ and the interior by
$I^\circ$. A
family
$\{F_p\}_{p\in I}$  of closed subsets in
$S^1$ is called an {\it intrinsic circle system on the interval $I$ } if it
satisfies the following
axioms:

\smallskip
\item{(I0)} The point $p$ is contained in $F_p$ for every $p\in I$.
\item{(I1)} If the set $F_p\cap F_q$ is non-empty, then $F_p=F_q$.
\item{(I2)} If $p'\in F_p$ and $q'\in F_q$
satisfy $p\prec q\prec p'\prec q'(\prec p)$, then $F_p=F_q$ holds.
\item{(I3)}
 Let $(p_n)_{n\in {\bf N}}$ and $(q_n)_{n\in {\bf N}}$ be two sequences in
$I^\circ$ such that
$\displaystyle\lim_{n\to \infty}p_n=p$
and $\displaystyle\lim_{n\to \infty}q_n=q$ respectively where $p,q\in I$.
Suppose that $q_n \in F_{p_n}$ for all $n$. Then
$q\in F_p$ holds.

\medskip
\noindent
{\it Remark.}
In [Um] intrinsic circle systems were defined on the whole circle $S^1$.
This new definition is a slight generalization to any  subinterval in $S^1$.
\medskip

We give two examples of intrinsic circle systems from the papers [Um] and
[TU 1].
\medskip

\ni{\bf Examples.} (i) Let $\gamma:S^1\to {\bf R}^2$ be a $C^2$ regular
simple closed  curve. Let $p$ be a point on $\gamma$ and
denote by $C^\bullet_p$ the largest circle
in the domain bounded by $\gamma$ that touches $\gamma$ in $p$. Set
$$
F^\bullet_p=\gamma\cap C^\bullet_p.
$$
It is easy to see that $\{F^\bullet_p\}_{p\in S^1}$ is an intrinsic circle
system. One can similarly define an intrinsic circle
system $\{F^\circ_p\}_{p\in S^1}$
using the smallest 
circle $C^\circ_p$ that is contained in the exterior domain of $\gamma$ and
touches
$\gamma$ in $p$ instead of
$C^\bullet_p$, see [Um]. If $\gamma$ is $C^3$, then it is easy to see that
the curvature of
$\gamma$ has a critical point at $p$ if either $F^\bullet_p$ or $F^\circ_p$
is connected.

(ii) Let $f:P^1\to P^1$ be a diffeomorphism of the
real projective
line. Let $p$ be a point in $P^1$ and denote by ${\cal P}_p$
the one-parameter
family of projective transformations of $P^1$ with the same
$1$-jet as $f$ in
$p$. We assume that ${\cal P}_p$ is parameterized by the real numbers and
consider $f\circ
P^{-1}_t$ for $P_t\in {\cal P}_p$. Then there are two numbers $t_0\le t_1$
such that 
 $f\circ P^{-1}_t$ has only a fixed point in $p$ if $t\not\in[t_0,t_1]$ and
more fixed points
than $p$ if $t\in(t_0,t_1)$. We assume the parameter to be chosen such that
$f\circ P^{-1}_t$
moves points locally on the left of $p$ away and brings those locally on the
right of $p$
closer if $t<t_0$. Let $F^\bullet_p$ denote the fixed point set of $f\circ
P^{-1}_{t_0}$ and
$F^\circ_p$ the fixed point set of $f\circ P^{-1}_{t_1}$. It is proved in
[TU 1] that
$\{F^\bullet_p\}_{p\in P^1}$ and $\{F^\circ_p\}_{p\in P^1}$ are
intrinsic circle systems. A point $p$ is called a {\it projective point of
$f$} if there is a
projective transformation that has the same $3$-jet as $f$ in $p$. In
general there is only a
projective transformation with the same $2$-jet as $f$ at a point $p$. If
either
$F^\bullet_p$ or $F^\circ_p$ is connected, then $p$ is a projective point.
\medskip

The following basic but easy lemma is proved in [Um] for an intrinsic circle
systems on $S^1$. The proof in the more general case is exactly the
same. Notice that the idea behind the lemma is essentially due to H. Kneser
[Kn], although not in
this abstract setting. One does not need axiom (I3) in the proof of the
lemma.

\medskip
\noindent
{\bf 3.1. Lemma. } {\it
Let $\{F_p\}_{p\in I}$ be an intrinsic circle system on
$I=[a,b]$.
Suppose that $F_a=F_b$ and $F_a \cap (a,b)$ is empty.
Then there exists a point $c\in (a,b)$ such that
$F_c$ is connected and contained in $(a,b)$.}\medskip

Lemma 3.1 applied to Example (i) has the classical Four-Vertex Theorem as a
consequence and is
nothing but a reformulation of its proof in [Kn]. Applied to Example (ii),
the theorem of Ghys
that a diffeomorphism of $P^1$ has at least four projective
points follows, see
[TU 1].
\bigskip

\noindent{\bf B. Intrinsic conic systems}\medskip

We will define an intrinsic conic system on an interval $I$ of $S^1$
to be a  set of functions
from $S^1\times S^1$ into the nonnegative even integers extended by $\infty$
that are indexed by a subset of $\bar I \times \bar I$ and
satisfy certain axioms. On one hand this is analogous to the
intrinsic circle systems defined above. On the other hand it is related to
divisors
on complex algebraic curves.
As we will see in Appendix C,
given a plane algebraic curve, one can consider the linear
system of divisors that come from intersections of the curve with
conics and use it to prove the formula of Cayley for the number of
sextactic points
mentioned in the introduction. The intrinsic conic systems that we consider
here do not correspond to the whole linear system, but only to those
coming from intersections with conics that are tangent to the curve
and do not cross it at any of the common points. This noncrossing property
is the reason why we restrict ourselves to even or infinite values of
the functions. See Example (i) that we give after the axioms and the
next section for full details of
this application. Generalizations to higher order intrinsic systems
and applications to Fourier series of periodic functions will be
given in [TU 2]. We explain a special case of the construction in [TU 2] in
Example (ii) after the axioms.

\medskip

Let $I$ be either the circle $S^1$ or an interval of $S^1$ that can be
open, closed or halfopen.
To avoid trivialities, we assume that the length of $I$ is positive.
 We set $$I^2_*:=(\bar I\times\bar I)\setminus \{(p,p)\;|\; p\in \bar
I\setminus I\},$$
i.e., $I^2_*$ is the closed square $\bar I \times \bar I$ with a corner
point $(p,p)$ removed if
$p\not \in I$. A  family $\{f_{p,q}\}_{(p,q)\in I^2_*}$
of functions $f_{p,q}:S^1\to 2{\bf N}_0\cup\{\infty\}$, where $2{\bf
N}_0$ denotes the nonnegative even integers,
is called an
 {\it intrinsic conic system
 on the interval $I$ } if it satisfies the axioms that will be listed
below.  (It is important
that the functions $f_{p,q}$ be defined on the whole circle $S^1$  since
axiom (A7) below might otherwise be violated in the applications.) Notice
that
$f_{p,q}$ is defined for
$p\ne q$ if and only if
$p,q\in
\bar I$ and
$f_{p,p}$ is defined if and only if $p\in I$. We
will denote the support of $f_{p,q}$ by $F_{p,q}$, i.e.,
$$F_{p,q}=\{r\in S^1\;|\; f_{p,q}(r)>0\}.$$ The value of
$f_{p,q}$ at a point $r$ will be called the {\it multiplicity of $r$ with
respect
to
$f_{p,q}$.} The sum over all values of $f_{p,q}$, which can of
course be infinite, is called the {\it total multiplicity of $f_{p,q}$}. A
point
$r$ in
$S^1$ will be called {\it sextactic} if its multiplicity with respect to
some $f_{p,q}$ is at
least six. We now list the axioms and follow them by examples  that
explain their  geometric
meaning.
\smallskip
\item{(A1)} $F_{p,q}$ is closed and $p,q\in F_{p,q}$ for all $(p,q)\in
I^2_*$.
\item{(A2)} $f_{p,q}=f_{q,p}$ for all $(p,q)\in I^2_*$.
\item{(A3)} If $F_{p,q}$ and $F_{p,r}$ have a point $s\ne p$ in common, then
$f_{p,q}=f_{p,r}$.
\item{(A4)}
If $p''\in F_{p,p'}$ and $q''\in F_{q,q'}$
satisfy $p\preceq q\preceq p'\preceq
 q'\prec p''\prec q''(\prec p)$ or
$p\succeq q\succeq p'\succeq
 q'\succ p''\succ q''(\succ p)$, and
$f_{p,p'}(p)\ge 4$ if $p=p'$, and
$f_{q,q'}(q)\ge 4$ if $q=q'$,
then $f_{p,p'}=f_{q,q'}$ holds.
\item{(A5)} If $f_{p,q}(r)\ge 4$ and $r\in I$, then $f_{r,r}=f_{p,q}$.
 \item{(A6)} Let $((p_n,q_n))$
be a sequence in $I^2_*$ such that
$\lim_{n\to \infty}(p_n,q_n)=(p,q)\in I^2_*$,
and let $(r_n^1)$ and $(r_n^2)$ be two sequences such that $r_n^i\in
\bar I\cap F_{p_n,q_n}$
and
$\lim_{n\to\infty}r_n^i=r\in \bar I$
for $i=1,2$. Assume $f_{p_n,q_n}(r_n^1)\ge k_1$ and $f_{p_n,q_n}(r_n^2)\ge
k_2$ for all $n$.
Then $f_{p,q}(r) \ge \max\{k_1,k_2\}$, and the inequality is strict if
$r_n^1$ and $r_n^2$
are different for all $n$.
\item{(A7)} The total multiplicity of $f_{p,q}$ is at least
six for all $(p,q)\in I^2_*$.\medskip
\item{(A8)}  
If $f_{p,q}(p)=2$, then $p$ is isolated in $F_{p,q}$.
\medskip

\ni{\bf Examples.} (i) Let $\gamma$ be a strictly convex curve in the affine
plan. We identify $\gamma$ with $S^1$. Let $C$ be a conic.
Then we can associate to $C$ a function on $S^1$ that
associates to a
point $r$ on $\gamma$ the multiplicity with which $C$ and $\gamma$ meet in
$r$. The multiplicity is
of course zero in points in which $C$ and $\gamma$ do not meet. Let $I$
denote $S^1$ or some interval
on $S^1$ and let $(p,q)\in I^2_*$. If
$p\ne q$, we let $C_{p,q}$ denote the maximal inscribed conic that is
tangent to
$\gamma$ in $p$ and $q$.   If $p=q$, we let
$C_{p,q}$ denote the maximal inscribed conic that meets $\gamma$ with
multiplicity at least four in
$p$. We let $f_{p,q}$ denote the function corresponding  to $C_{p,q}$ as
explained above. We
will prove in section three that $\{f_{p,q}\}_{(p,q)\in I^2_*}$ is an
intrinsic conic system.
The sextactic points of $\{f_{p,q}\}_{(p,q)\in I^2_*}$ are precisely the
globally maximal
sextactic points of $\gamma$.

(ii) For a real valued $C^4$-function $u$ on $S^1$, $n\ge
0$, one defines the
{\it osculating polynomial
$\varphi_s$} ({\it of order five}) {\it at a point $s\in S^1$} to be the
unique trigonometric
polynomial of degree two,
$$
\varphi_s(t)=a_0+a_1\cos t+b_1\sin t+a_2\cos 2t+b_2\sin 2t,$$ whose
value and first four derivatives at $s$ coincide with those of $u$ at $s$.
If $u$ is $C^5$ and the value and the first five derivatives of
$u$ and $\varphi_s$ coincide in $s$, i.e., if $\varphi_s$ hyperosculates $u$
in $s$, then we call $s$ a
{\it flex of $u$} ({\it of order five}). The existence of six flexes of
order
five can easily be proved as a consequence of the well-known fact that a
function has at least six zeros if its Fourier coefficients
$a_i$ and $b_i$ vanish for $i\le 2$, see [TU 2]. One can use intrinsic conic
systems
to prove the much stronger result that there are six
such flexes satisfying the {\it global property} that the osculating
polynomials $\varphi_s$
in the flexes 
support  $u$, i.e., either $\varphi_s\le u$ or $u\le \varphi_s$. The
intrinsic conic 
systems are defined as follows. To simplify the definition we assume $u$ to
be $C^\infty$.
Let $(p,q)\in S^1\times S^1$. If $p\ne q$, we
let $\varphi_{p,q}$ denote the smallest trigonometric polynomial of degree
two that is greater
or equal to $u$ and has the same values as $u$ in $p$ and $q$.  If $p=q$ we
let $\varphi_{p,q}$
denote the smallest trigonometric polynomial of degree two that is greater
or equal to $u$ and
has the same
$1$-jet as $u$ in $p=q$.
We now define $f_{p,q}:S^1\to 2{\bf
N}_0\cup\{\infty\}$ by setting $f_{p,q}(r)=2k$ if the $2k-1$-jets of $u$ and
$\varphi_{p,q}$
agree in $r$ but not the $2k$-jets, $f_{p,q}(r)=0$ if the values of $u$ and
$\varphi_{p,q}$ do
not agree in $r$, and $f_{p,q}(r)=\infty$ otherwise. One can now prove that
$\{f_{p,q}\}_{(p,q)\in S^1\times S^1}$ is an intrinsic conic system and
similarly define an
intrinsic conic system using trigonometric polynomials that are smaller or
equal to $u$. We
refer to [TU 2] for much more general results concerning Fourier polynomials
of arbitrary
degree. Notice that Fourier polynomials of degree one lead to intrinsic
circle system.

\medskip
We now start deriving consequences of the axioms of an intrinsic conic
system.
The following lemma is trivial.\medskip

\ni{\bf 3.2. Lemma.} {\it
Let $\{f_{p,q}\}_{(p,q)\in I^2_*}$ be an intrinsic conic system on
an interval $I$ in $S^1$.
Then for any  subinterval $J$ of $I$,
the restriction $\{f_{p,q}\}_{(p,q)\in J^2_*}$ is an
intrinsic conic system on $J$.
}\qed

The following lemma is an immediate consequence of axiom (A7).\medskip

\ni{\bf 3.3. Lemma.} {\it If $F_{p,p}$ only consists of the point $p$, then
$p$ is
sextactic.}\qed

We also have the following two easy lemmas.\medskip

\ni{\bf 3.4. Lemma.} {\it  $f_{p,p}(p)\ge 4$ for every $p\in I$.}\medskip

\ni{\it Proof.} Let $(p_n)$ and $(q_n)$ be two sequences in $I$ that
converge to $p$ and $p_n\ne q_n$
for all $n$. Applying (A6) to the situation $r_n^1=p_n$ and $r_n^2=q_n$, we
get $f_{p,p}(p)>2$
since $p_n\in F_{p_n,q_n}$ and $q_n\in F_{p_n,q_n}$ by (A1). Using that
$f_{p,p}(p)$ is an
even number, we get $f_{p,p}(p)\ge 4$.
\qed

\ni{\bf 3.5. Lemma.} {\it If $(p_n)$ and $(q_n)$ are sequences in $I^\circ$
that both
converge to
$p\in I$,
$q_n\in F_{p_n,p_n}$, and $q_n$ is different from $p_n$,
then $p$ is sextactic. }\medskip

\ni{\it Proof.} By Lemma 3.4, we have $f_{p_n,p_n}(p_n)\ge 4$ for every
$n$. Hence
(A6) implies that $f_{p,p}(p)> 4$ since we can choose $r^1_n=p_n$ and
$r^2_n=q_n$.
It follows that $f_{p,p}(p)\ge 6$ and hence that $p$ is sextactic.\qed

\ni{\bf 3.6. Lemma.} {\it If $r\in F_{p,q}\cap I$ is not isolated in
$F_{p,q}\cap I$, then
$r$ is a sextactic point with infinite multiplicity with respect to
$f_{p,q}$. }\medskip

\ni{\it Proof.} We assume that $f_{p,q}(r)$ is a finite number $k$.
 Let $(r_n)$
be a sequence in
$F_{p,q}$ of pairwise different points that are all different from
$r$ and converges to $r$. We now apply (A6) to the situation $p_n=p$,
$q_n=q$, $r_n^1=r$ and
$r_n^2=r_n$. It follows that $f_{p,q}(r)>k$, a contradiction. Hence
$f_{p,q}(r)=\infty.$ \qed

If the support of $f_{p,p}$ is a connected set, then we say that $p$ is a
{\it clean
sextactic point}. Lemmas 3.3 and  3.6 imply that a clean sextactic point in
$I$
is sextactic and
that moreover every point in the intersection of the support of $f_{p,p}$
with
$I$ is a clean sextactic point if $p$
is.\medskip

\ni{\bf 3.7. Lemma.} {\it If $F_{p,q}\ne F_{q,q}$, then $f_{p,q}(q)=2$, $q$
is
isolated in $F_{p,q}$ and $F_{p,q}$ has at least two components.}
\medskip

\ni{\it Proof.} {If $F_{p,q}\ne F_{q,q}$, then (A5) implies that
$f_{p,q}(q)=2$.  By (A8),
we know that $q$ is isolated in $F_{p,q}$. Since $q$ must of course be
different from $p$, we
see that $F_{p,q}$ must have at least two connected components.} \qed

We set
$$
F^*_{p,q}=\left\{
\matrix{F_{p,q}\hfill &\qquad \hbox{if} \,\,\,f_{p,q}(p)\ge 4,\cr
F_{p,q}\setminus\{p\} &\qquad \hbox{if} \,\,\,f_{p,q}(p)= 2.\cr
}
\right.
$$

\medskip
\noindent
{\bf 3.8. Lemma.} {\it
Let $I$ be $S^1$ or an interval on $S^1$ and
$\{f_{p,q}\}_{(p,q)\in I^2_*}$ an intrinsic conic system
on $I$. Then for each $p\in I$,
 $\{F^*_{p,q}\}_{q\in\bar I}$
is an intrinsic circle system on $\bar I$.}

\medskip
\noindent
{\it Proof.} First notice that $f_{p,q}$ is defined for all $q\in\bar I$
since $p\in I$.
To prove (I0), notice that $q\in F_{p,q}^*$ if $q\ne p$. If $q=p$,
then $f_{p,q}(p)=f_{p,p}(p)\ge 4$ and hence $q\in F_{p,q}=F_{p,q}^*$.

To prove (I1) assume that $s\in F^*_{p,q}=F^*_{p,r}$. If $s=p$, then
$f_{p,q}(p)\ge 4$ and $f_{p,r}(p)\ge 4$ and hence $f_{p,q}=f_{p,s}$ by
(A5). Now assume that $s\ne p$. Then (A3) implies that  $f_{p,q}=f_{p,s}$.

Property (I2) follows from (A4).

To prove that (I3) holds, let $(q_n)$ and $(r_n)$ be sequences in $I^\circ$
with limits $q$ and $r\in \bar I$  respectively, and assume that
$r_n\in F_{p,q_n}^*$, then it follows from (A6) that $r\in F_{p,q}$ and
(I3) follows if $r\ne p$. If
$r=p$, then we have to prove that
$f_{p,q}(p)\ge 4$. If only finitely many of the $r_n$ are equal to $p$,
then (A6)
applied to $r_n$ and the constant sequence $p$ implies that $f_{p,q}(p)\ge
4$. If infinitely many of the
$r_n$ are equal to $p$, then for these $r_n$ we have $f_{p,q_n}(p)\ge 4$
and hence that $f_{p,q}(p)\ge 4$. \qed
\medskip

\noindent
{\bf 3.9. Lemma.} {\it
Let $I=[a,b]$ or $[a,b)$, $a\prec b$, be an interval on $S^1$,
$\{f_{p,q}\}_{(p,q)\in I^2_*}$  an intrinsic conic system
on $I$ such that $f_{a,a}=f_{a,b}$ and $F_{a,b}\cap (a,b)$ is empty.
Then there exists a point $c\in (a,b)$ such that
$f_{a,c}=f_{c,c}$ and $F_{c,c}$ is contained in $[a,b)$. Furthermore, $a$ is
isolated in $F_{a,c}$ and $F_{a,c}$ has exactly two components.

Similarly, if $I=[a,b]$ or $(a,b]$,  $f_{b,b}=f_{a,b}$
and $F_{a,b}\cap (a,b)$ is empty, then there exists a point $c\in (a,b)$
such that
$f_{b,c}=f_{c,c}$ and $F_{c,c}$ is contained in $(a,b]$. Furthermore, $b$ is
isolated in $F_{b,c}$ and $F_{b,c}$ has exactly two components. }
\medskip

\noindent
{\it Proof.} We only prove the first part; the second part is similar.
By Lemma 3.8 we know that  $\{F^*_{a,q}\}_{q\in \bar I}$ is an intrinsic
circle system
on $\bar I$. Clearly it satisfies $F^*_{a,a}=F^*_{a,b}$.
Hence there exists a point $c\in (a,b)$ by Lemma 3.1 such that
$F^*_{a,c}$ is connected and contained in $(a,b)$. Hence $a\not\in
F^*_{a,c}$
which implies that $F_{a,a}\ne F_{a,c}$ and $F_{a,c}$ has exactly two
connected components.
Therefore  we have $F_{a,c}=F_{c,c}$ and hence $f_{a,c}=f_{c,c}$. It is
clear
that
$F_{c,c}$ is contained in
$[a,b)$.
\qed

\medskip

\medskip
\noindent
{\bf 3.10. Lemma.} {\it
Let $I=[a,b]$ or $[a,b)$, $a\prec b$, be a closed interval of $S^1$ and
$\{f_{p,q}\}_{(p,q)\in I^2_*}$ be an intrinsic conic system
on $I$. Suppose that $f_{a,a}=f_{a,b}$ and $F_{a,b}\cap (a,b)$ is empty.
Then there exist two distinct points $a_1\prec b_1$
in $(a,b)$
such that $f_{a_1,b_1}=f_{a_1,a_1}$ and $F_{a_1,b_1}= F_{a_1,a_1}\subset
(a,b)$, $b_1$ is isolated in
$F_{a_1,b_1}$
 and $F_{a_1,b_1}$
has exactly two  components.

Similarly, if  $I=[a,b]$ or $(a,b]$, $f_{b,b}=f_{a,b}$ and $F_{a,b}\cap
(a,b)$ is
empty, then there exist two distinct points $a_1\prec b_1$
in $(a,b)$
such that $f_{a_1,b_1}=f_{b_1,b_1}$ and $F_{a_1,b_1}= F_{b_1,b_1}\subset
(a,b)$, $a_1$ is isolated in
$F_{a_1,b_1}$
 and $F_{a_1,b_1}$
has exactly two  components. }

\medskip
\noindent
{\it Proof.} We prove the first part of the lemma.
By Lemma 3.9, there exists a point $b_1\in (a,b)$ such that
$f_{a,b_1}=f_{b_1,b_1}$ with support in $I$ and $a$ is isolated in
$F_{a,b_1}$.
We can assume that  $b_1$ is such that $(a,b_1)\cap F_{a,b_1}=\emptyset$.
Then by Lemma 3.9,
there exists a point $a_1\in (a,b_1)$ such that
$f_{a_1,b_1}=f_{a_1,a_1}$ with support in $(a,b_1]\subset (a,b)$.
\qed

The next two propositions will be the main tools to find sextactic points
in sections four
and five.\medskip

\ni{\bf 3.11. Proposition.}  {\it
Let $I$ be a closed or halfopen interval of $S^1$ with endpoints $a$ and
$b$  and let
$\{f_{p,q}\}_{(p,q)\in I^2_*}$ be an intrinsic conic system
on $I$. Suppose that
 $F_{a,b}\cap (a,b)$ is empty. We assume furthermore that either
$f_{a,b}=f_{a,a}$ or
$f_{a,b}=f_{b,b}$ holds (at least one of these conditions makes sense when
$I$ is
halfopen).
Then there is a sextactic point
$r$ in
$(a,b)$.}\medskip

\ni{\it Proof.} Let $J=[a_1,b_1]$ be an interval as in Lemma 3.10.
Let ${\cal C}_J$ denote the set of $(\alpha,\beta)\in (a,b)\times (a,b)$
such that $\alpha\ne \beta$,
$f_{\alpha,\beta}=f_{\alpha,\alpha}$ with support in $J$,
$F_{\alpha,\beta}\cap
I(\alpha,\beta)=\emptyset    $, and $F_{\alpha,\beta}$ consists of
precisely two components, one of
which is the isolated point
$\beta$. Here $I(\alpha,\beta)$ denotes the open interval with endpoints
$\alpha$ and $\beta$.
(Notice that we do not assume that $\alpha<\beta$.) We know from Lemma 3.10
that
${\cal C}_J$ is nonempty.

We let $\delta_{\alpha,\beta}$ denote  the distance between $\alpha$ and
$\beta$. Let $\delta$
denote the infimum over
$\delta_{\alpha,\beta}$ for
$(\alpha,
\beta)\in {\cal C}_J$.

We consider
a sequence
$\{(\alpha_n,\beta_n)\}$ in ${\cal C}_J$ such that
$\delta_{\alpha_n,\beta_n}$ converges to $\delta$.
By going to subsequences if necessary, we may assume that
$$
\lim_{n\to\infty}\alpha_n=\alpha,\qquad
\lim_{n\to\infty}\beta_n= \beta. $$
If $\alpha=\beta$, then it follows immediately from Lemma 3.5 that $\alpha$
is a sextactic point.
We can therefore assume that $\delta>0$.
To simplify the notation we will assume that $\alpha\prec \beta$.
By (A6), we have $f_{\alpha,\alpha}(\alpha)\ge 4$ and
$f_{\alpha,\beta}(\alpha)\ge 4$ and hence
$f_{\alpha,\alpha}=f_{\alpha,\beta}$ by (A5). (We do not claim that
$F_{\alpha,\beta}$ is contained in
$J$.) We can assume that $\alpha$ and $\beta$ are isolated in
$F_{\alpha,\beta}$ since otherwise we
have a sextactic point by Lemma 3.6.
Let $\beta'$ be the point in $F_{\alpha,\beta}\cap
(\alpha,\beta]$ closest to $\alpha$.  We now apply Lemma 3.9 to the
interval $[\alpha,\beta']$ and we
find a point $\gamma\in (\alpha,\beta')$ such that $(\gamma, \alpha)\in
{\cal C}_J$. Clearly
$\delta_{\gamma,\alpha}<\delta$, which is a contradiction. Hence there is a
sextactic point in $J\subset (a,b)$.
\qed

\medskip
\noindent{\bf 3.12. Proposition.} {\it Let
$\{f_{p,q}\}_{(p,q)\in I^2_*}$ be an intrinsic conic system
on $I$, where $I$ is some interval on $S^1$ or the whole circle $S^1$.
Assume that $p$ and $q$
are contained in distinct components of $F_{p,q}$ and that there is a third
component of $F_{p,q}$
between $p$ and $q$ on $I$.
Then there is a point
$r\in I$ such that $r\ne q$ and $F_{r,q}$  has two connected components one
of which is $\{q\}$.
}

\medskip
\noindent {\it Proof.} We can assume without loss of generality that there
is a point
$p'$ in $F_{p,q}\cap I$ different
from both
$p$ and $q$  such that the open interval $J$ between $p$ and
$p'$ on $I$ does not meet $F_{p,q}$.
 By Lemma 3.8, $\{F^*_{q,x}\}_{x\in I}$ is an intrinsic
circle system on $I$.
Since $p, p'\in F^*_{q,p}$ there is by Lemma 3.1
 a point $r$ in $J$ such that $F^*_{q,r}$ is
connected and contained in $J$. Then $F_{q,r}$ has two connected
components one of which is $\{q\}$. \qed
\medskip

The last two propositions were the main technical results of this section.
We use them to prove the following theorem.\medskip

\ni{\bf 3.13. Theorem.} {\it Let
$\{f_{p,q}\}_{(p,q)\in S^1\times S^1}$ be an intrinsic conic system
on $S^1$. Then $\{f_{p,q}\}_{(p,q)\in S^1\times S^1}$ has at least three
sextactic points.}\medskip

\ni{\it Remark.} This theorem is optimal as the intrinsic conic system
$\{f^\bullet_{p,q}\}_{(p,q)\in S^1\times S^1}$ of a contractible branch of
a real regular cubic shows, see section four and Appendix C.\medskip

\ni{\it Proof.} We first prove the existence of one sextactic point.
Let $p$ be a point on $\gamma$
that we can assume not to be sextactic. Then $f_{p,p}(p)=4$.
 Hence $p$ is isolated in $F_{p,p}$. We
therefore have a point $q$ in $F_{p,p}$ that is different from $p$
and such that the open interval
$(q,p)$ does not meet $F_{p,p}$. Proposition 3.11 now implies that there is
a sextactic
point $s$ in the interval
$(q,p)$.

To prove that there are two further sextactic points we proceed as follows.
Let $r$ be some point
different from $s$.  If $F_{r,s}$ consists of two components we
have two sextactic points different from $s$ by Proposition 3.11. If
$F_{r,s}$ consists of  three
components at least, we can use Proposition 3.12 to find a point $r'$ such
that $F_{r',s}$ consists of two components,
and the existence of the two new sextactic points follows again from
Proposition 3.11.
\qed

\bigskip
\noindent
{\bf \S 4\hskip 0.1in An application to strictly convex curves.}\medskip

In this section we use the theory of intrinsic conic systems to give a
complete proof of the results of Mukhopadhyaya in [Mu 2] on the existence
of inscribed and
circumscribed osculating conics of strictly convex curves in an affine
plane.  As was
pointed out in section two, such a curve is the same thing as a simple
closed curve in $P^2$
without inflection points.\medskip

\ni{\bf 4.1. Theorem} (Mukhopadhyaya). {\it Let $\gamma$ be a strictly
convex curve in the affine
plane $A^2$. Then
$\gamma$ has at least three circumscribed osculating conics and at least
three inscribed
osculating ellipses. In particular, $\gamma$ has at least six sextactic
points that are
globally maximal or minimal.}\medskip

The osculating conic at a point $p$ of a curve $\gamma$ is an ellipse if
and only if the affine
curvature of $\gamma$ at $p$ is positive. It therefore follows from the
theorem that a strictly
convex curve must have points with positive affine curvature.\medskip

We first point out that the theorem has an interesting corollary which does
not seem to follow from the
other proofs of the existence of sextactic points. We will denote the open
disk which a Jordan
curve $\gamma$ in $A^2$ bounds by $D_\gamma$ and refer to it as the {\it
interior domain of
$\gamma$}. We let $\kappa_M$ denote the maximum of the affine curvature of
$\gamma$, $\kappa_m$ its
minimum and $A(D_\gamma)$ the area of $D_\gamma$.\medskip

\ni{\bf 4.2. Corollary.} {\it Let $\gamma$ be a strictly convex curve in
the affine plane $A^2$.
Then $\kappa_M>0$ and
$$\pi \kappa_M^{-3/2}\le A(D_\gamma)$$
with an equality if and only if $\gamma$ is an ellipse. If the affine
curvature of $\gamma$ is
positive then we also have

$$A(D_\gamma)\le \pi \kappa_m^{-3/2}$$
with an equality if an and only if $\gamma$ is an ellipse.
} \medskip

\ni{\it Remark.} Both inequalities follow from exercises 4 and 15 in
section \S 27 of [Bl 2] if the
affine curvature is positive.\medskip

\ni{\it Proof of the Corollary.} We already observed that $\kappa_M>0$. Now
let $C$ be one of the
inscribed osculating ellipses and denote its affine curvature by $\kappa$.
Then $\kappa\le
\kappa_M$. The area of the interior domain of $C$ is
$\pi\kappa^{-3/2}$. The first inequality follows immediately. The second
inequality follows by
arguing similarly with one of the circumscribed conics which must be an
ellipse since the affine curvature is
positive.\qed

Before proving  Theorem 4.1 we need to introduce the relevant intrinsic
conic systems.\medskip

 We say that the ellipse $C_1$ is  contained in the ellipse $C_2$ if
$D_{C_1}\subset D_{C_2}$ where $D_{C_i}$ is the interior domain of $C_i$.
An inscribed ellipse is
said to be {\it maximal} if it is not
strictly contained in any other inscribed ellipse. \medskip

Let $p,q$ be two different points on $\gamma$. Let $\Gamma_{p,q}$ be the
one dimensional
family of ellipses  that is tangential to $\gamma$ in $p$ and $q$.  In one
direction, this
family converges to the  closed line segment
$\overline{pq}$.  Since $\overline{pq}$ meets $\gamma$ transversally in $p$
and
$q$, we have  inscribed ellipses in the family. Thus there is a unique
maximal inscribed ellipse
 in the family
$\Gamma_{p,q}$ that we will denote by $C^\bullet_{p,q}$. \medskip

 We can also define the maximal inscribed ellipse
$C^\bullet_{p,q}$ when $p=q$.
Fix a point $p$ on
$\gamma$.
Let $\Gamma_{p,p}$ be the one dimensional family of ellipses
that is tangential to $\gamma$ in $p$ with multiplicity at least four. The
osculating ellipse
of $\gamma$ in $p$ is defined since $p$ is not an inflection point and it
is of course
contained in the family
$\Gamma_{p,p}$. No ellipse in
$\Gamma_{p,p}$ can cross
$\gamma$ in
$p$ except possibly the osculating ellipse. In one direction, this family
$\Gamma_{p,p}$ converges
to the point $p$. In that same direction after passing the osculating
ellipse, the ellipses lie
locally around $p$ inside of $\gamma$.  Hence we have inscribed ellipses in
the family.
There is therefore a unique maximal inscribed
ellipse
 in the family $\Gamma_{p,p}$ that we denote by $C^\bullet_{p,p}$.\medskip

For any pair of points $(p,q)\in S^1\times S^1=\gamma\times\gamma$, let
$f^\bullet_{p,q}:S^1\to 2{\bf N}_0\cup\{\infty\}$ denote the function that
associates to a point
$r\in S^1 =
\gamma$ the multiplicity with which $C^\bullet_{p,q}$ and $\gamma$ meet in
$r$. If $r\not\in
C^\bullet_{p,q}$, then of course $f^\bullet_{p,q}(r)=0$.\medskip

The following  is obvious:\medskip

\ni  {\it  The functions $f^\bullet_{p,q}$  satisfy axioms (A1), (A2),
(A3) (A5)  and  (A8)  for
intrinsic conic systems for all $(p, q)\in S^1\times S^1$.}\qed\medskip

Notice that axiom (A4) is an easy consequence of the fact that two
ellipses cannot
meet in more than four points without being identical. Hence we
have:\medskip

\ni  {\it  The functions $f^\bullet_{p,q}$  satisfy axiom (A4) for
intrinsic conic systems
for every pair
$(p,q)\in S^1\times S^1$.} \qed\medskip

Assume that $C$ is an ellipse that meets $\gamma$ in a point $p$ with
multiplicity two. Then
$C$ and $\gamma$ do not cross in $p$ and there is another ellipse $C^\prime$
tangent
to $\gamma$ in $p$ and containing
$C$ which locally around $p$ lies between $\gamma$ and $C$. This implies
the following
lemma.
\medskip

\noindent {\bf  4.3. Lemma.} {\it Let $p$ and $q$ be two distinct points on
$\gamma$. Then the
maximal ellipse $C^\bullet_{p,q}$ can be characterized as the only inscribed
ellipse that
meets $\gamma$ in $p$ and $q$ and satisfies
 one of the following two properties.
\item{(i)} $C^\bullet_{p,q}$ meets
$\gamma$  in at least three different points.
\item{(ii)}
$C^\bullet_{p,q}$ meets $\gamma$ with multiplicity
at least four
either in $p$ or in $q$ and it does not have any other points in common with
$\gamma$.

\ni In particular, $C^\bullet_{p,q}$ meets $\gamma$ with
total multiplicity
six at least.
}\qed
\medskip

The following lemma follows from the fact that if $C$ is an ellipse in
$\Gamma_{p,p}$ that
lies locally around $p$ inside of $\gamma$ and meets $\gamma$ in $p$
precisely with
multiplicity four, then there is a different ellipse $C^\prime$ in
$\Gamma_{p,p}$ that contains
$C$,  is contained in the osculating ellipse at $p$, lies locally around $p$
inside of $\gamma$
and also meets $\gamma$ in $p$ with  precisely multiplicity  four.

\medskip
\noindent
{\bf 4.4. Lemma.} {\it Let $p$ be an arbitrary point on $\gamma$. Then the
maximal ellipse
$C^\bullet_{p,p}$ can be characterized as the only inscribed ellipse that
meets $\gamma$ in $p$
with multiplicity at least four and satisfies one of the following two
properties.
\item{(i)} $C^\bullet_{p,p}$ meets $\gamma$   in at least two different
points.
\item{(ii)} $C^\bullet_{p,p}$ meets $\gamma$ in $p$ with multiplicity at
least six and has no other
points in common with $\gamma$. It follows that $C^\bullet_{p,p}$ is the
osculating ellipse at
$p$ and $p$ is a sextactic point.

\ni In particular, $C^\bullet_{p,p}$ meets $\gamma$ with
total multiplicity
six at least. }\qed
\medskip

We have the following immediate corollary of Lemmas 4.3 and 4.4. \medskip

\ni  {\it  The functions $f^\bullet_{p,q}$  satisfy axiom (A7) for
intrinsic conic systems
for every pair
$(p,q)\in S^1\times S^1$. }\qed\medskip

It is therefore only left to prove that axiom (A6) is satisfied.\medskip

\ni{\bf 4.5. Lemma.} {\it The functions $f^\bullet_{p,q}$  satisfy axiom
(A6) for
intrinsic conic systems for every pair
$(p,q)\in S^1\times S^1$.}\medskip

\ni {\it Proof.} We consider a small interval around $r=\gamma(t_0)$ on the
curve $\gamma$ and
parallel coordinates $(x,y)$ in which $\gamma$ on this interval corresponds
to points on the
$x$-axis. Let $r^1_n=\gamma(t_n^1)$ and $r^2_n=\gamma(t_n^2)$ be two
sequences converging to $r$
and assume that $r^1_n$ and $r^1_n\in C^\bullet_{p_n,q_n}$ where $p_n$ and
$q_n$ converge to $p$
and $q$ respectively. We assume that $f^\bullet_{p_n,q_n}(r^i_n)\ge k_i$
for $i=1,2$.
Assume that $k_1\le k_2$.
We can write $C^\bullet_{p_n,q_n}$ locally around $r$ as a graph of a
function $g_n(t)$ in the
parallel coordinates for sufficiently big $n$. We can assume after going to
a subsequence if
necessary that the $C^\bullet_{p_n,q_n}$ converge to an ellipse $C$ that is
the graph of a function
$g$ in the parallel coordinates.
We have that $g_n(t_n^1)=g_n'(t_n^1)=\dots
=g_n^{(k_1-1)}(t_n^1)=0$ and  $g_n(t_n^2)=g_n'(t_n^2)=\dots
=g_n^{(k_2-1)}(t_n^2)=0$. By taking
limits it clearly follows that
 $g(t_0)=g'(t_0)=\dots =g^{(k_2-1)}(t_0)=0$. Hence $C$ and $\gamma$ meet in
$r$ with multiplicity $k_2$ at
least.

We now prove that $C$ and $\gamma$ meet in $r$ with multiplicity greater
than $k_2$
  when $r_n^1$ and $r^2_n$ are different for all $n$. Set
$i=k_2-k_1$. There is by Rolle's Theorem of Elementary Calculus for every
$n$ a $s_n^1$ between
$t_n^1$ and $t^2_n$ such that $g_n^{(k_1)}(s_n^1)=0$. Similarly we find an
$s^2_2$ between $s_n^1$
and
$t_n^2$ with $g_n^{(k_1+1)}(s_n^2)=0$, and inductively an $s^j_n$ between
$s^{j-1}_n$ and $t_n^2$
with $g_n^{(k_1+j-1)}(s_n^j)=0$ for $j=1,\dots, i+1$. Taking limits we
obtain
 $g(t_0)=g'(t_0)=\dots =g^{(k_2)}(t_0)=0$. This proves that $C$ and
$\gamma$ meet in $r$
with multiplicity greater than
$k_2$.

The claim in the lemma now follows after we prove that $C^\bullet_{p,q}$
meets $\gamma$
in $r$ at least with the same multiplicity as $C$. Notice that $C$ contains
the points
$p$ and $q$ and is inscribed in $\gamma$. If $p\ne q$, then it follows that
$C^\bullet_{p,q}$
lies between $\gamma$ and $C$ since it is maximal with this property. Hence
$C^\bullet_{p,q}$
meets $\gamma$ at least with the same multiplicity in $r$ as $C$ if $p\ne
q$. If $p=q$, the same follows
if we can show that $C$ meets $\gamma$ with multiplicity at least four in
$p=q$. If there are
infinitely many $n$ such that $p_n=q_n$, then this follows as in the first
paragraph of the proof.
If there are infinitely many $n$ such that $p_n\ne q_n$, this follows as in
the second
paragraph of the proof.
\qed

With
help of circumscribed conics, we next associate in an analogous manner a
second intrinsic conic
system to a strictly convex curve $\gamma$.  For this purpose it will be
more convenient to
assume that we are in $P^2$, since otherwise we would for example need to
take both branches of a
hyperbola into account when defining an interior domain. A conic $C$ in
$P^2$ bounds a closed disk
$D_C$ and we say that it circumscribes a simple closed curve $\gamma$ if
$\gamma\subset D_\gamma$.
It is clear what we mean by a minimal circumscribed conic.
\medskip

Let $p,q$ be two different points on $\gamma$. Let $\Gamma_{p,q}$ be
the one dimensional family of conics  that is tangential to $\gamma$ in $p$
and $q$.  In one
direction, this
family converges to the union of the tangent lines of $\gamma$ at $p$ and
$q$.  Since
$\gamma$ is strictly convex, we have circumscribed conics in the family.
(Working in an affine
plane, we might not have a circumscribed ellipse in this family. This
happens e.g.~in points where the affine
curvature is nonpositive.) Thus there is a unique minimal circumscribed
conic in the family
$\Gamma_{p,q}$ that we will denote by
$C^\circ_{p,q}$.
\medskip

 We now define $C^\circ_{p,q}$ in the case that $p$ and $q$ coincide.
For $p$ on
$\gamma$ we  let $\Gamma_{p,p}$ be the one dimensional family of conics
that is tangential to $\gamma$ in $p$ with multiplicity at least four. We
have circumscribed conics in the family since $\gamma$ is strictly convex.
There is therefore a
unique minimal circumscribed conic in the family $\Gamma_{p,p}$ that we
denote by
$C^\circ_{p,p}$.\medskip

Now for any pair of points $(p,q)\in S^1\times S^1=\gamma\times\gamma$, we
define
$f^\circ_{p,q}:S^1\to 2{\bf N}_0\cup\{\infty\}$ to be the function that
associates to a point
$r\in S^1 =
\gamma$ the multiplicity with which $C^\circ_{p,q}$ and $\gamma$ meet in
$r$. \medskip

We have already seen that $\{f_{p,q}^\bullet\}$ is an intrinsic conic
system. The same
arguments imply that $\{f_{p,q}^\circ\}$ is an intrinsic conic system. Thus
we have the
following proposition.

\medskip
\ni{\bf 4.6. Proposition.} {\it Both $\{f_{p,q}^\bullet\}$ and
$\{f_{p,q}^\circ\}$
are intrinsic conic systems on $S^1$.}\qed\medskip

We can now give a proof of Theorem 4.1. \medskip

\ni{\bf 4.7. Proof of Theorem 4.1.} The claim is an immediate consequence of
the last
proposition and Theorem 3.13, as well as the definition of the relevant
intrinsic conic systems
with help of inscribed and circumscribed conics. \qed

We remind the reader that a sextactic point $p$ on a simple closed curve
$\gamma$ is
called a {\it clean
maximal (resp.~minimal) sextactic point} if
 $C_{p,p}^\bullet\cap \gamma(S^1)$ (resp.~$C_{p,p}^\circ\cap \gamma(S^1)$)
is
connected. The following theorem will be proved in [TU 2].
\medskip

\medskip
\noindent
{\bf 4.8. Theorem.} {\it
Let $\gamma:S^1\to  A^2$ be a strictly convex curve.
Then $\gamma$ has at least three clean maximal and at least three clean
minimal
sextactic points.}

\medskip 
In the following we will prove Lemma 2.3 only assuming
$C^4$-differentiability
of the arc $\gamma$. This together with Remark (ii) after Proposition 5.1
should
make it clear how to prove the theorems in the introduction for
$C^4$-curves.

\medskip
\noindent
{\bf 4.9.  Proof of Lemma 2.3.}  We will assume here that
$\gamma:[0,1]\to P^2$ is a {\it simple} $C^4$-arc and explain after the
proof of
Lemma 4.10 how this follows from the other assumptions that we make on
$\gamma$.

We assume that $\gamma$ is free of inflection
points.  Instead of assuming that $\gamma$
is free of sextactic points, we make the following weaker assumptions: We
assume
that the osculating conic $C_t$ at $\gamma(t)$ crosses $\gamma$ in
$\gamma(t)$
for every $t\in(0,1)$. We also assume that $\gamma$ enters the interior
domain of
$C_0$ in $\gamma(0)$ and that it lies locally outside of $C_1$ in
$\gamma(1)$.
The claim of the lemma follows if we can show that
$C_0$ does not meet
$C_1$. We therefore assume that $C_0$ and $C_1$ meet.  If $\gamma|_{(0,1]}$
meets $C_0$, then we let $c\in(0,1]$ be the smallest number such that
$\gamma(c)\in C_0$, and we set $\hat\gamma=\gamma|_{[0,c]}$. If
$\gamma|_{(0,1]}$ does not meet $C_0$, then we extend the arc $\gamma$ by
continuing on $C_1$ up to the first point where $C_1$ meets $C_0$. We denote
the extended arc by $\hat\gamma$ and assume it to be parameterized on the
interval $[0,c]$. We can assume that $\hat\gamma$ is a $C^4$-regular curve.
Now we set
$$
\sigma=\hat\gamma\cup C_0|_{[\hat\gamma(0),\hat\gamma(c)]}.
$$
and assume that $\sigma$ is parameterized on $[0,1]$ with
$\sigma(0)=\sigma(1)=\hat\gamma(c)$. Notice that $\sigma$ is simple and can
be
assumed to be a $C^4$-regular arc that makes a loop which possibly does
not close smoothly in $\sigma(0)=\sigma(1)$. Notice that the interior
angle in $\sigma(0)=\sigma(1)$ is less than or equal to $\pi$.  Such a curve
$\sigma$ cannot exist because of the following lemma.
\qed

\medskip
\noindent
{\bf 4.10 Lemma.} {\it Let $\sigma:[0,1]\to A^2$ be a strictly convex simple
closed curve in the affine plane $A^2$ that is $C^4$-regular everywhere in
$[0,1]$, but possibly so that $\dot\sigma(0)\not=\dot\sigma(1)$. Denote by
$D$ the closed domain bounded by $\sigma$. Set $p=\sigma(0)=\sigma(1)$ and
suppose that the interior angle $\theta$ at $p$ is less than or equal to
$\pi$.
Suppose moreover that the affine curvature function $\kappa$ of $\sigma$
is non-decreasing and non-constant. Then there exists a sextactic
point $s\in (0,1)$ such that the osculating conic at $s$
does not not coincide with the osculating conics at $\sigma(0)$
and $\sigma(1)$. }
\medskip

\noindent
{\it Proof.}\hskip .3cm First notice that the osculating conic $C_0$ at
$\gamma(0)$ is not inscribed in $D$ since we are assuming that the curve
$\sigma$ has non-decresing and non-constant affince curvature functon
and hence there exists a point
$\varepsilon\ge 0$ such that the closed arc
$\sigma([0,\varepsilon])$ (or point if $\varepsilon=0$)
is a connected component of $C_0\cap
\sigma$
and $\sigma$ enters
the interior domain $D$ of $C_0$ in $\sigma(\varepsilon)$.

 For each $s,t\in (0,1)$,
let $C^\bullet_{s,t}$ be the maximal conic inscribed in $D$ with
$\sigma(s),\sigma(t)\in C^\bullet_{s,t}$ and $C^\bullet_{s,t}$ meeting
$\sigma$ in $\sigma(s)$ with multiplicity four if $s=t$.  Now for any pair
of points $(s,t)\in (0,1)\times (0,1)$, we
define $f^\bullet_{s,t}: S^1\to 2{\bf
N}_0\cup\{\infty\}$ as follows:
 We set $f^\bullet_{s,t}(p)=2$ if
 $p$ is on $C^\bullet_{s,t}$, otherwise we set
$f^\bullet_{s,t}(p)=0$. (Notice that  $C^\bullet_{s,t}$ and
$\sigma$ might meet in $p$ if $\theta=0$.)
For  $r\in S^1\setminus \{p\} =\partial D\setminus \{p\}$,
$f^\bullet_{s,t}$ is the
multiplicity with which $C^\bullet_{s,t}$ and $\sigma$ meet in $r$ if it is
less than five, otherwise we set $f^\bullet_{s,t}(r)=\infty$.
Then $f^\bullet_{s,t}$ satisfies the axioms of an intrinsic conic system  on
every closed interval $[a,b]$ such that $0<a<b<1$.

Assume that  $\theta<\pi$.
We fix two distinct points $t_0,s_0 \in (0,1)$.
We set $C=C^\bullet_{t_0,s_0}$.
Then $C$ meets $\sigma|_{(0,1)}$ with total multiplicity six.
Applying Proposition 3.12, we find a point $u\in (0,1)$
such that $F_{t_0,u}$ consists of two
components one of 
which is $t_0$.
By Proposition 3.11, we find a sextactic point $s$
between $t_0$ and $u$ whose osculating
conic $C_s$ is inscribed.
Since the interior angle is less than
$\pi$ by assumption, $C_s$ cannot pass through $p$.
This implies $C_s\ne C_1,C_0$.

 Next we consider the case $\theta=\pi$.
We set
$$
\delta={\rm inf}\{t\in (0,1]\,\;\,\sigma(t)\in C_1 \}.
$$
If $\delta=0$, there is a sequence $(u_n)_{n\in N}$
converging to zero such that $\sigma(u_n)\in C_1$,
and hence $C_0=C_1$, a contradiction.
So $\delta>0$.

We fix two distinct points $t_0,s_0 \in (0,1)$ such that
$t_0<s_0$ and $t_0\in (0,\delta)$.
We set $C=C^\bullet_{t_0,s_0}$.

First we consider the case that
$C$ meets $\sigma|_{(0,1)}$ with total multiplicity six.
We set
$$
t_0^*={\rm inf}\{t\in (0,1]\,\;\,\sigma(t)\in C=C^\bullet_{t_0,s_0} \}.
$$

Applying Proposition 3.8 and Lemma 3.1 to the
 intrinsic circle
system
$((F^\bullet)^*_{t_0^*,p})_{p\in S^1}$, we find a point $u\in (0,1)$
such that $C_{t_0^*,u}=C^\bullet_{u,u}$.
By Proposition 3.11, we find a sextactic point $s$ between $t_0^*$
and $u$ whose osculating conic $C_s$ is inscribed.
Since $C_0$ is not inscribed in $D$, $C_s\ne C_0$ is obvious.
If $C_1$ is not inscribed, $C_s\ne C_1$ also holds.
So we may assume $C_1$ is inscribed.
Suppose $\kappa(u)=\kappa(1)$. Then  $C^\bullet_{u,u}=C_u=C_1$ holds,
where $C_u$ is the osculating conic at $u$.
Since $C_{t_0^*,u}=C^\bullet_{u,u}$, we have $C_{t_0^*,u}=C_1$,
which contradics
the fact that $C_1$ does not pass through $\sigma(t_0^*)$.
So we have $\kappa(u)<\kappa(1)$ and hence
$\kappa(s)<\kappa(1)$, which implies that $C_s\ne C_1$.

Finally, we consider the case that
 $C$ meets $\sigma|_{(0,1)}$
only at $t_0$ and $s_0$ and that the multiplicity is equal to two
in both points. Then $C$ must also be tangent to $\sigma$ at the point $p$.
We fix two points $t_1\in (0,t_0),\,\,s_1\in (t_0,s_0)$, and
set $C'=C^\bullet_{t_1,s_1}$.
If $C'$ passes through $p$, then $C'$ meets $C$ with
total multiplicity five at least and we have $C'=C$, which is a
contradiction.
So $C'$ does not pass through $p$.
Thus $C'$ meets $\sigma|_{(0,1)}$ with total multiplicity six.We can now apply the arguments in the last
paragraph to $C'$ instead of $C$ to find a sextactic point $s$ whose
osculating conic is
inscribed and different from $C_0$ and $C_1$.
\qed

\bigskip
\noindent
{\bf \S 5\hskip 0.1in An application to simple closed curves.}\medskip

In this section we prove the theorems in the introduction except Part (i)
of Theorem 1.2
which is Mukhopadhyaya's Theorem that we already proved in the last section.
Theorem 1.1 is
the same as 5.2, Theorem 1.2 (ii) is in 5.3 and 5.5, (iii) is in 5.4,
Theorem 1.3 is in 5.4 and 5.5.\medskip

We  start with a proposition that will be our main technical tool.

\medskip
\noindent
{\bf 5.1. Proposition.}  {\it
Let $\sigma:[0,1]\to P^2$ be an arc of a curve $\hat\sigma: [-\epsilon,
1+\epsilon]\to P^2$ for
$\epsilon>0$ that does not have any self-intersections. We assume that
$\sigma(t)$ is not an inflection point for any $t\in (0,1)$
and that $\sigma(0)$ and $\sigma(1)$ are either inflection or
minimal sextactic points of $\hat\sigma$.
Then there exists a sextactic point on $\sigma|(0,1)$. }\medskip

\noindent{\it Remark.} (i) One sees from the proof below that there is a
maximal
sextactic point on $\sigma|(0,1)$ if both $\sigma(0)$ and $\sigma(1)$ are
inflection points.

(ii) In the proof of Case (b) in the proof below we will be dealing with a
curve that is only $C^5$
at one point and otherwise smooth. One can see directly that axiom (A6) is
satisfied at this point.
Notice that one can define an intrinsic conic system for a strictly convex
curves $\gamma$ that
is only $C^4$ as follows. In the notation of section four one sets
$$
f^\bullet_ {p,q}(r)=
\left\{\matrix{
2 & {\rm if}\; C^\bullet_{p,q}\; {\rm meets }\; \gamma\;{\rm at }\; r\;
{\rm with\; multiplicity
}\; 2,\hfill\cr
4 & {\rm if}\; C^\bullet_{p,q}\; {\rm meets }\; \gamma\;{\rm at }\; r\;
{\rm with\; multiplicity
}\; 4,\hfill\cr
\infty\; & {\rm if}\; C^\bullet_{p,q}\; {\rm meets }\; \gamma\;{\rm at }\;
r\; {\rm with\;
multiplicity\; higher\; than\; four}.\cr }\right.
$$
Then $\{f^\bullet_ {p,q}(r)\}$ satisfies the axioms of an intrinsic conic
system.

(iii)  We explain here how one can easily prove a weak version of
Proposition 5.1 under
generic assumptions on the arc $\sigma$ using affine curvature.
 Assume 
that
$\sigma:[0,1]\to A^2$ is a regular arc with no inflection points in
$\sigma(0,1)$ and that the
endpoints
$\sigma(0)$ and $\sigma(1)$ are inflection points with the property that the
tangent lines
there only meet $\sigma$ with finite multiplicity. Then we will show below
that the open arc
$\sigma(0,1)$ contains a sextactic point.
Fabricius-Bjerre [Fa] makes this observation under the stronger assumption
that $\sigma$ meets
the tangent lines in the endpoints with multiplicity three precisely and
uses it to
prove a weak version of Theorem 1.1.

To prove the claim in the previous paragraph, we choose
coordinates $(x,y)$ in $A^2$
such that $\sigma(0)$ corresponds to $(0,0)$
and the $x$-axis is the oriented tangent line
of $\sigma$ at $t=0$.  After
reparameterizing $\sigma$ we can write it as a
graph $y=y(x)$ for $x\ge 0$.
Since $\sigma$ is of finite type, we have
that
$$y(x)=\alpha x^n+o(x^n).$$ We can assume that $\alpha>0$
by either changing the orientation of $\sigma$ or the sign of the
$y$-coordinate.
 The
affine
curvature $\mu(x)$ can be expressed as
$$
\mu(x)=-{1\over
9(y^{''})^{8/3}}(5(y^{'''})^2-3y^{''''}y^{''}),
$$
 see [Bl 2], p.14,
formula (83). A short calculation shows that
$$
\lim_{t\to 0+}\mu(t)=-\infty.
$$ 
Similarly it follows that
$$\lim_{t\to
1-}\mu(t)=-\infty.
$$
As a consequence there is a point on $\sigma(0,1)$ where the affine
curvature 
takes on its maximum value and this point is then the sextactic point whose
existence we wanted
to show.

\medskip
\noindent
{\it Proof of Proposition 5.1.} Assume that both $\sigma(0)$ and $\sigma(1)$
are minimal
sextactic points. If there is
no sextactic point on $\sigma|(0,1)$, then by Lemma 2.3 all osculating
conics along $\sigma|(0,1)$
are disjoint and it follows that the osculating conic at $\sigma(0)$ must
contain $\sigma|[0,1]$ in
its interior domain since $\sigma(0)$ is a minimal sextactic point. In
particular, the osculating
conic at $\sigma(1)$ is contained in the interior domain of the osculating
conic at $\sigma(0)$. We
can reverse the roles of $\sigma(0)$ and $\sigma(1)$ in this argument and
prove that the osculating
conic at $\sigma(0)$ is contained in the interior domain of the osculating
conic at $\sigma(1)$
which is a contradiction. Hence the proposition is proved if both
$\sigma(0)$ and $\sigma(1)$ are
minimal sextactic points. We will therefore assume in the rest of the proof
that at least one of
the points $\sigma(0)$ or $\sigma(1)$ is an inflection point.

In the rest of the proof we will denote the tangent line of  $\hat\sigma$
in $\sigma(0)$ by $L_0$ and
the one in
$\sigma_1(1)$ by
$L_1$. We will assume $L_0$ and $L_1$ parameterized such that the tangents
of $L_0$ and $\sigma$
coincide in
$\sigma(0)$ as well as those of $L_1$ and $\sigma$ in $\sigma(1)$.

The following three cases can occur, see Figure 5.1.

\medskip
\item{(a)} The curve $\sigma((0,1))$ neither meets  $L_0$
nor  $L_1$.
\item{(b)}
The curve $\sigma((0,1))$  intersects $L_0$.
\item{(c)}
The curve $\sigma((0,1))$  intersects $L_1$.
\medskip
\centerline{\epsfxsize=3in \epsfbox{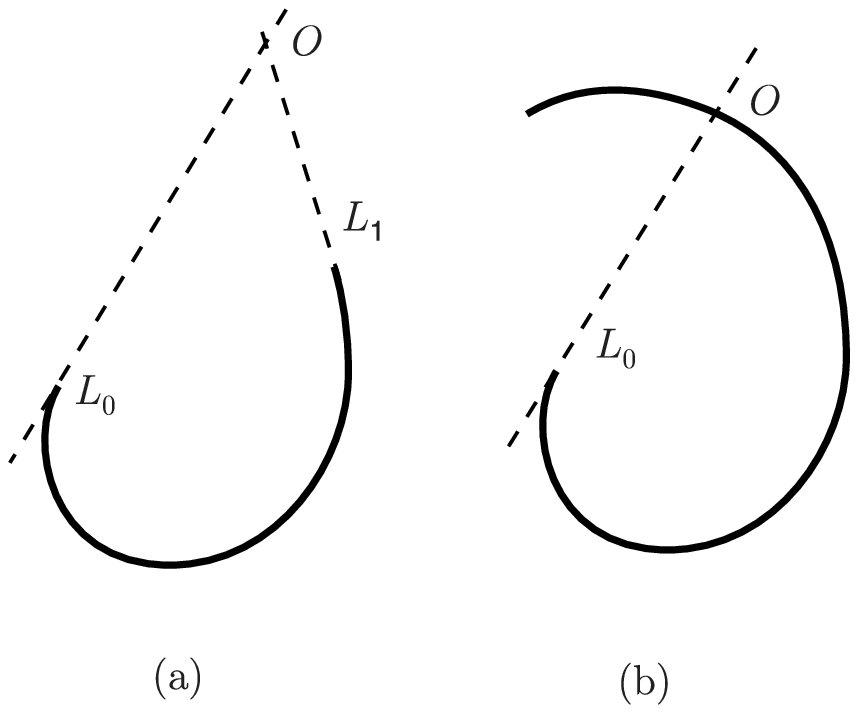}}
\centerline{Figure 5.1}
\bigskip

Notice that case (b) and case (c) are up to orientation of the curve
identical.  \medskip

\smallskip
\noindent
({\it Case }(a))
The tangent line $L_0$ at $\sigma(0)$ and
the tangent line $L_1$ at $\sigma(1)$ meet in a point that we denote by $O$.
(See Figure 5.1.)

Consider the simple closed curve
$$
\gamma:=L_0|_{[O,\sigma(0)]}\cup \sigma([0,1])\cup
L_1|_{[\sigma(1),O]}.
$$
We will first prove that the curve $\gamma$ is nullhomotopic.  It follows
from Theorem 2.2
that there is an
 affine plane $A^2$ that contains $\sigma([0,1])$ (but not necessarily
$\gamma$).
Let $L$ be the line segment in this affine plane between $\sigma(0)$ and
$\sigma(1)$.
The triangle $L_0|_{[O,\sigma(0)]}\cup L\cup L_1|_{[\sigma(1),O]}$ is
contractible.
Hence $\gamma$ is homotopic to $\sigma([0,1])\cup L$ which is nullhomotopic
since it is contained in
an affine plane.

Let $D_\gamma$ denote the closed disk bounded by $\gamma$. Notice that
$L_0$ and $L_1$ do not meet any
interior point of $D_\gamma$. One can move $L_0$ (or $L_1$) slightly so
that it does not meet
$D_\gamma$. It follows that $\gamma$ and $D_\gamma$ lie in an affine plane
$A^2$ and that
 $D_\gamma$ is convex.

 We first assume that both $\sigma(0)$ and $\sigma(1)$ are inflection
points.
As in section four,  we  consider inscribed conics.
Let $(p,q)$ be a pair of different points on  $\sigma([0,1])$.
We let  $C_{p,q}^\bullet$ denote the maximal inscribed conic
that lies in $D_\gamma$ and is tangential
to $\gamma$ in $p$ and $q$. If $p\in \sigma((0,1))$, then $p$ is not an
inflection point of $\sigma$
and we can define $C_{p,p}^\bullet$ as the maximal inscribed conic that
meets  $\sigma$ with
multiplicity at least four in
$p$. We define $f^\bullet_{p,q}(r)$ for $r\in \sigma$ to be the
multiplicity with which
$\hat\sigma$ and
$C^\bullet_{p,q}$ meet. If $r$ lies on the open segments between $O$ and
either $\sigma(0)$ or
$\sigma(1)$, then we set $f^\bullet_{p,q}(r)=2$ if $r\in C_{p,q}^\bullet$,
otherwise we set
$f_{p,q}(r)=0$.  One can prove exactly as in the last section that
$\{f_{p,q}\}$ is an intrinsic conic system on the open interval
$\sigma((0,1))$.  In fact (A8) follows from the fact that
 $C^\bullet_{\sigma(0),p}$ (resp. $C^\bullet_{p,\sigma(1)}$)
never touches $L_0\setminus\{\sigma(0)\}$
(resp. $L_1\setminus\{\sigma(1)\}$).

We will now show that
 $C^\bullet_{\sigma(0),\sigma(1)}$ meets $\sigma((0,1))$ in
a point $r$. After having shown this the claim of the proposition follows
from 3.12 and 3.11 in the
case we are now considering.

Assume that there is no such
point
$r$, i.e.,
$C_{\sigma(0),\sigma(1)}^\bullet$  only meets $\gamma$ in $\sigma(0)$ and
$\sigma(1)$. Since
$\sigma(0)$ and $\sigma(1)$ are inflection points we have that
$C_{\sigma(0),\sigma(1)}^\bullet$ and
$\hat\sigma$ can only meet with multiplicity two both in $\sigma(0)$ and
$\sigma(1)$. Of course
$C_{\sigma(0),\sigma(1)}^\bullet$ also meets $L_0$ and $L_1$ with
multiplicity two in $\sigma(0)$
and $\sigma(1)$. Hence we can increase the conic
$C_{\sigma(0),\sigma(1)}^\bullet$ in such a way that locally around
$\sigma(0)$ and $\sigma(1)$ it
stays inside of $D_\gamma$. Since $C_{\sigma(0),\sigma(1)}^\bullet$ only
meets $\gamma$ in
$\sigma(0)$ and $\sigma(1)$ this is not only true locally around
$\sigma(0)$ and $\sigma(1)$, but
globally, contradicting the maximality of
$C_{\sigma(0),\sigma(1)}^\bullet$. Hence $\gamma$ and
$C_{\sigma(0),\sigma(1)}^\bullet$ meet in a third point $r$ which must lie
on $\sigma((0,1))$
since $C_{\sigma(0),\sigma(1)}^\bullet$ can not meet
the tangent lines $L_0$ and $L_1$ except in $\sigma(0)$ and $\sigma(1)$.
This finishes the proof
when both $\sigma(0)$ and $\sigma(1)$ are inflection points.

Now assume that $\sigma(1)$ is minimal sextactic. We can assume that
$\sigma(0)$ is an inflection point as pointed out at the beginning of the
proof.  We assume that there is no
sextactic point on
$\sigma((0,1))$. This implies that $\sigma((0,1))$ lies in the interior
domain of the osculating
conic at
$\sigma(1)$ that we will denote by $C$. Notice that $L_1$ is tangent to $C$
in $\sigma(1)$, but
does not meet it otherwise. The conic $C$ enters $D_\gamma$ in $\sigma(1)$
and leaves it in a point
$O'$ on $L_0$ that lies between $\sigma(0)$ and $O$.
The curve $\sigma([0,1])\cup C|_{[\sigma(1),O']}$ is $C^\infty$-regular
except in $\sigma(1)$
where it is $C^5$.  Furthermore it satisfies the condition in case (b). We
will prove below that
 the closed curve
$$
\hat\gamma:=L_0|_{[O',\sigma(0)]}\cup \sigma([0,1])\cup
C|_{[\sigma(1),O']}
$$
has an inscribed osculating conic. Such a conic can only be osculating at
points in $\sigma((0,1))$. It now
follows that we have a maximal sextactic point on
$\sigma((0,1))$ contradicting our
assumption. Notice that it does not follow in this case that we have a
maximal sextactic
point since we only prove its existence assuming that there is no
sextactic points in $\sigma((0,1))$.

If $\sigma(0)$ is minimal sextactic, we can of course use the same argument
thus finishing the proof
of Case (a).
\smallskip
\noindent
({\it Case }(b)) We can assume that only one of $\sigma(0)$ and $\sigma(1)$
is
minimal sextactic as observed at the beginning of the proof. Assume that
$\sigma(0)$
is minimal sextactic and that there is no
sextactic point on
$\sigma((0,1))$.   Then $\sigma([0,1])$ would lie in the closed interior of
the
osculating conic at
$\sigma(0)$ and
$L_0$ could not meet $\sigma((0,1))$ which contradicts that we are in case
(b).
We therefore have a sextactic point on $\sigma((0,1))$ if $\sigma(0)$ is
minimal
sextactic.
Hence we can assume that $\sigma(0)$ is an inflection point.
 Let
$O$ be the point where
$L_0$ meets
$\sigma$ for the first time. (See Figure 5.1.) Consider the simple closed
curve
$$
\gamma:=L_0|_{[O,\sigma(0)]}\cup \sigma|_{[\sigma(0),O]}.
$$
In the following we need to include the possibility
that $\sigma|_{[\sigma(0),O]}$ is the $C^5$-curve we met in case (a).
Notice that $\gamma$ bounds a closed convex domain that we denote by
$D_\gamma$.
Exactly as in Case (a), we define for a pair of points
 $(p,q)$ on  $\sigma|_{[\sigma(0),O]}$ such that
$(p,q)\ne (\sigma(0),\sigma(0)), p\ne O$, and $q\ne O$,
the maximal inscribed conic $C_{p,q}^\bullet$. We also define
$f^\bullet_{p,q}$ as in Case (a)
and prove that $\{f^\bullet_{p,q}\}$ is an intrinsic conic system on
any halfopen arc $(\sigma(0),p]\subset (\sigma(0),O)$ of $\gamma$.
(Here we have
to check axiom (A6) separately for the point
where the curve from case (a) is only $C^5$.)
Fix an arbitrary point $p$ on $(\sigma(0),O)$.
Consider the conic $C_{\sigma(0),p}^\bullet$.
Since $\sigma(0)$ is an inflection point and the angle at $O$ is acute,
$C_{\sigma(0),p}^\bullet$ must by arguments as in the proof of
Case (a) either meet a point $r(\ne p)$ on $\sigma$ between
$\sigma(0)$ and $O$ or it meets $\sigma$ in $p$ with multiplicity at least
four.
Now it follows from Propositions 3.12 and 3.11 that there is a
maximal sextactic point on the open arc
of
$\sigma$  between $\sigma(0)$ and $O$. This
proves Case (b).
\smallskip
\noindent
({\it Case }(c)) Same proof as in Case (b).
\qed

\medskip
\ni{\bf 5.2. Theorem.} {\it Let $\gamma:S^1\to P^2$ be a simple closed
curve without
self-intersections that is not nullhomotopic. Then $\gamma$ has at least
three sextactic
points.}\medskip

\ni{\it Proof.} By  Theorem 2.1, due to M\"obius, there are at least three
intervals of true
inflection points on $\gamma$. We therefore find
 three different arcs on $\gamma$ whose endpoints are inflection points.
Now the claim of
the theorem follows immediately from Proposition 5.1. \qed

\medskip
\noindent
{\bf 5.3. Proposition.} \it
Let $\gamma:S^1\to P^2$ be a simple closed curve that is nullhomotopic
and meets every line in $P^2$.
Then $\gamma$ has at least four sextactic points.
\rm

\medskip
\noindent
{\it Proof.} By Theorem 2.2, the curve $\gamma$ must
have at least four intervals of true inflection points, and hence at least
four different subarcs
whose endpoints are inflection points. The claim now follows from
Proposition
5.1. \qed

\medskip
\noindent
{\bf 5.4. Theorem.} {\it
Let $\gamma:S^1\to A^2$ be a regular closed convex curve.
Then it has at least two sextactic points. The total number of inflection
and
sextactic points is at least four. In particular, a convex curve with one
inflection point has at
least three sextactic points.}

\medskip
\noindent
{\it Proof.} We proved in section three that $\gamma$ has at least six
sextactic points if it has
no inflection points, i.e., if $\gamma$ is strictly convex.

If there is more than one interval of inflection points,
then we find at least two sextactic points by Proposition 4.1. We have
therefore proved the theorem except when the set of inflection points on
$\gamma$
is an interval.

Fix two distinct non-inflection  points $a$ and $b$ on $\gamma$.
Consider the minimal circumscribed conic $C_{a,b}^\circ$ that touches
$\gamma$ in $a$ and $b$.
If $C_{p,q}^\circ$ touches $\gamma$ with multiplicity four in either $a$ or
$b$, we consider the
intrinsic conic system $\{f_{p,q}^\circ\}$ on the closed interval between
$a$ and $b$
which does not contain an inflection point. Then by Proposition 3.11, we
find a minimal sextactic
point between
$a$ and
$b$. If $C_{a,b}^\circ$ meets $\gamma$ both in $a$ and $b$ with
multiplicity two, then there must
be a third point $c$ on $\gamma$ lying in $C_{a,b}^\circ$. This point
cannot be an inflection
point since
$\gamma$ lies inside of $C_{a,b}^\circ$ around $c$.   After renaming the
points $a$, $b$ and $c$ we
can assume that $c$ lies on the open arc between $a$ and $b$ that is free
of inflection points. We now
consider the intrinsic conic system $\{f^\circ_{p,q}\}$ on the closed
interval $[a,b]$ and
use Propositions 3.12 and 3.11 to prove the existence of a minimal
sextactic point between $a$ and
$b$. In both cases we have a minimal sextactic point and a connected set of
inflection points on the
curve
$\gamma$. There are therefore two open intervals on the curve that are
bounded by a minimal
sextactic point and an inflection point. Now it follows from Proposition
5.1 that we have a
sextactic point on each interval. We have hence proved in this case that
there are at least three
sextactic points on
$\gamma$ and one inflection point. This finishes the proof of the theorem.
\qed

\ni{\bf 5.5. Theorem.} {\it
Let $\gamma:S^1\to A^2$ be a regular closed curve that is not convex.
Then it has at least three sextactic points. The total number of inflection
and sextactic points is
at least six. In particular, if $\gamma$ is not convex and has two
inflection points, then it has
at least four sextactic points.}

\medskip
\ni{\it Proof.} The curve $\gamma$ not being convex has true inflection
points. If
$\gamma$ has at least three intervals of inflection points, then it has at
least three sextactic
points by Proposition 5.1.  So we may assume that
$\gamma$ has exactly two  intervals of inflection points that divide
$\gamma$ into two arcs.
Let $\sigma$ denote the boundary of the convex hull of $\gamma$.
Since we are assuming that $\gamma$ has exactly two intervals inflection
points,
$\sigma$ consists of an arc of $\gamma$ and a line segment between
two points
$a$ and $b$ on $\gamma$. We choose two arbitrary different points $p$ and
$q$ 
on the open arc of $\gamma$ between
$a$ and $b$ and
  consider the  minimal circumscribed conic $C^\circ_{p,q}$ touching
$\sigma$  in $p$ and $q$.
Then the conic $C^\circ_{p,q}$ will either meet $\gamma$
in $p$ or $q$ with multiplicity four or it will meet
$\gamma$ in a third point $r$.
When the second case occurs,
by replacing $p$ by $r$ if necessary,
we may assume that $r$ lies between $p$ and $q$.
 In both cases minimal circumscribed conics touching
$\sigma$ along $\gamma$
gives rise to an intrinsic
conic system on a
closed interval of
$\gamma$ between $p$ and $q$. Then  we can deduce the existence of a
minimal sextactic point on
the arc between $a$ and $b$ using Propositions 3.11 and 3.12.
 Proposition
5.1 now implies that there
are two further sextactic points on the arcs of $\gamma$ between the
inflection points and the
minimal sextactic point. There is a fourth sextactic point on the arc
between the inflection points that
lies inside the convex hull of $\sigma$. This finishes the proof of
the theorem.
\qed

\vfill
\eject

\ni{\bf Appendix A: Simple closed curves with few inflection points}\medskip

We will use our paper [TU 1] to prove the following theorem which is the
second part of
Theorem 2.2.\medskip

\ni{\bf A.1. Theorem.} {\it
A contractible simple closed curve $\gamma: S^1\to P^2$ with less than four
intervals of true
inflection points is contained in an affine plane.}\medskip

The first part of Theorem 2.2 is a consequence since a regular simple arc
$\sigma: [0,1]\to P^2$
such that  no $\sigma(t)$ for $t\in (0,1)$ is an inflection point can be
shown to be a part of a
simple closed curve $\gamma$ with at most two inflection points.
\medskip

Theorem A.1 will immediately follow from a result on curves on $S^2$. In
fact, the preimage of
$\gamma$ under the canonical projection $p:S^2\to P^2$ consists of two
simple closed curves since
$\gamma$ is contractible. Let $\hat\gamma$ be one of these curves. Then the
other one is
$-\hat\gamma$.
\medskip

A point on $\hat\gamma$ is called an {\it inflection point}
if the osculating circle at that point is a great circle. The inflection
points
on $\gamma$ and $\hat\gamma$
correspond since $p$ maps great circles to lines. We can define a {\it true
inflection point} of
$\gamma$ as a point where the geodesic curvature of $\hat\gamma$ changes
sign. An interval on
$\hat\gamma$ is said to consists of true inflection points if the geodesic
curvature vanishes
and changes sign there.  Also the true inflection points on $\gamma$ and
$\hat\gamma$ correspond.

\medskip
Theorem A.1 now follows immediately from the following theorem after
observing that a lift
$\hat\gamma$ of the curve $\gamma$ cannot contain a great semicircle as a
subarc.

\medskip
\noindent
{\bf A.2. Theorem.} {\it
Let $\hat\gamma$ be a simple closed curve on $S^2$ with at most
three intervals of true inflection points. Then $\hat\gamma$ lies in a
closed hemisphere.
Moreover, $\hat\gamma$ lies in an open hemisphere if and
only if it does not contain any great semicircle as a subarc.}

\medskip
\noindent
{\it Proof.} The claim that $\hat\gamma$ lies in a closed hemisphere was
already proved
in section two of [TU 1], see also the arguments later in this proof.
It is therefore sufficient to show that
$\hat\gamma$ lies in an open hemisphere if and
only if it does not contain any
great semicircle as a subarc.

Assume that  $\hat\gamma$ contains a great semicircle $J$ as a subarc.
Then any great circle on $S^2$ meets $J$ and hence $\hat\gamma$ too.
So $\hat\gamma$ can not lie any open hemisphere.

We now prove the converse.
Suppose that $\hat\gamma$ is a simple closed curve with
at most three intervals of true inflection points and does not contain
any  great semicircle as a subarc.
By [TU 1], there exist  four points
$t_1<t_2<t_3<t_4$ on $\hat\gamma$ such that $t_1,t_3$ are clean maximal
vertices and $t_2,t_4$ are clean minimal vertices.
The osculating planes of $\hat\gamma$, considered as a space curve, at
these four points bound a
simplex $S$ in ${\bf R}^3$ containing $\hat\gamma$, see [TU 1]. If the
origin
of ${\bf R}^3$ lies in
the interior of this simplex, the curve $\hat\gamma$ has at least four
intervals of true inflection
points, see [TU 1], which contradicts our assumption.
Thus, the origin lies in the boundary of the simplex $S$.
Hence the origin lies in an osculating plane $P_j$ at one of the
clean vertices $t_j$.
 We set $C_j=P_j\cap S^2$. Then $C_j$ is a great circle which is
an osculating circle of $\hat\gamma$ at $t_j$. The curve $\hat\gamma$ lies
completely
on one side of $P_j$ in ${\bf R}^3$ and hence also on one side
of $C_j$ on $S^2$, which was the proof
in [TU 1] that $\hat\gamma$ lies in a closed hemisphere.
 Since $t_j$ is a clean vertex we have that
$C_j\cap \hat\gamma$ is connected.  Now we use that $\hat\gamma$ does not
contain a great semicircle as
a subarc. It follows that $C_j\cap \hat\gamma$
 is a point or a great circle arc
whose length is less than $\pi$.
Now let $C_j^+=C_j|_{(S,N)}$ be an open  great semicircle  bounded by
two points $S,N$ on $C_j$ such that $C_j \cap \hat\gamma$ is contained in
$C_j^+$. Rotate  the great circle $C_j$ slightly around the axis in ${\bf
R}^3$
passing through $S$ and $N$ away from $\hat\gamma$ into a great circle
$\hat C$. If the
rotation is sufficiently small, the curve
$\hat\gamma$ does not meet $\hat C$. Hence $\hat\gamma$ lies in an open
hemisphere, and we have
finished the proof.
\qed

\bigskip

\ni{\bf Appendix B: Examples of curves with few sextactic points}
\medskip

In this section we give examples of simple closed curves in the affine
plane with few sextactic
points that show together with the next appendix that the theorems in the
introduction are
optimal. Two of these examples are due to  Izumiya and Sano [IS] who came
up with them
in their  study of affine evolutes of convex curves.
\medskip

Part (ii) of Theorem 1.2 is optimal by Example B.1. Part (iii) of Theorem
1.2 is optimal by
Example B.2. Theorem 1.3 is optimal by Examples B.2 and B.4.
\medskip

We know from the proofs in section five that a simple closed curve with
more than three
(intervals) of inflection points has more than three sextactic points. If
it has three inflection
points, then it has at least three sextactic point. This is optimal by
Example
B.1. If it is convex and has
two inflection points, it has at least two sextactic points. This is optimal
by
Example B.2. If it is
not convex and has two inflection points, we know from Theorem 5.5 that it
has  at least four sextactic
points. This is optimal by Example B.3. If it is convex and has one
inflection point, then it has at least
three sextactic points by Theorem 5.4. This is optimal by Example B.4.
\medskip

\ni{\bf B.1 Example.} Here we give an example of a simple closed curve in
the affine plane that is
not convex, has three inflection points (two of which are true inflection
points) and only three
sextactic points. This shows that part (ii) of Theorem 1.2 is optimal.

We identify the affine plane with the complex plane and consider the map
$$z(t)=t+{3i\over 1+t^2} \qquad {\rm for}\qquad t\in {\bf R}.$$
We then get our example by setting
$$\gamma(t)=1/z(t) \qquad {\rm for }\qquad t\in {\bf R}\cup\{\infty\}.$$
Notice that the curve is regular in $t=\infty$. In Cartesian coordinates
the curve can be
expressed as
$$x(t)={t\over r(t)}\ ,\qquad y(t)={-3\over r(t)(1+t^2)}$$
where
$$r(t)=t^2+{9\over (1+t^2)^2}.$$
There is a sketch of the curve in Figure
B.1, where the inflection points and the sextactic points
are marked by $I$ and $S$ respectively.
\medskip

\def\A{\vtop{\halign{\hfill##\hfill\cr
               {\epsfxsize=3in \epsfbox{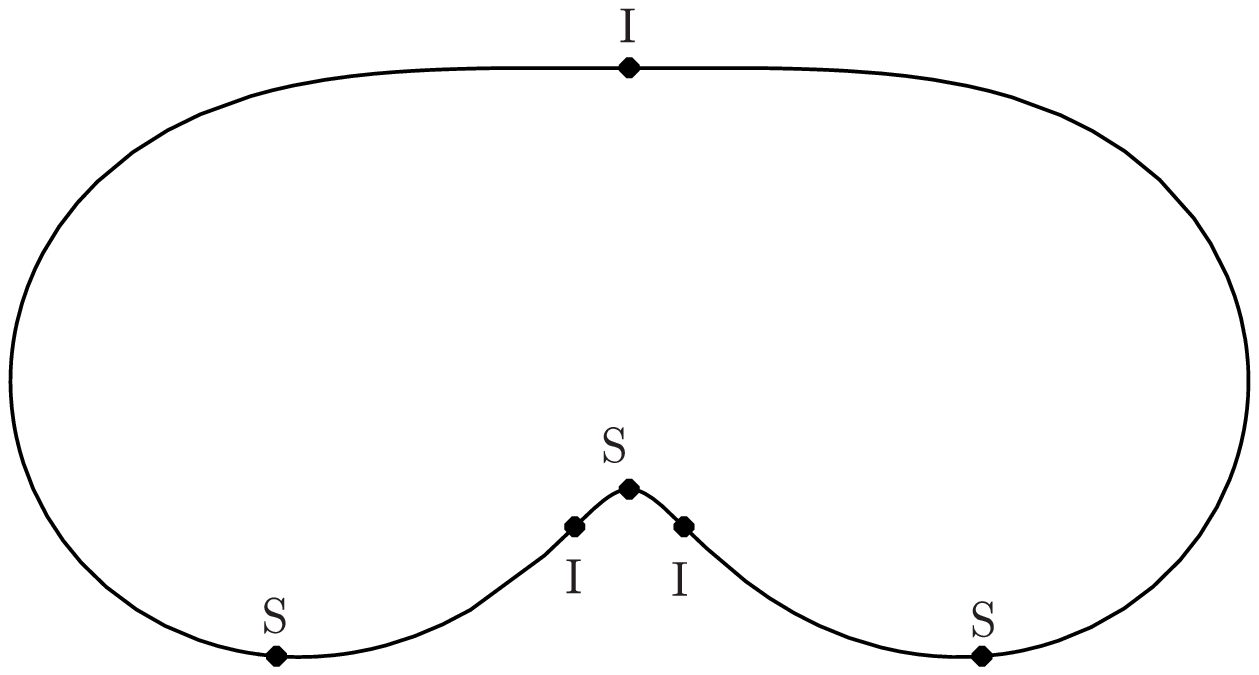}} \cr
            \noalign{\smallskip}{Figure B.1} \cr}}}
\centerline{\A}
\bigskip
\bigskip

\ni{\bf B.2 Example.} The following example is due to Izumiya and Sano
[IS].  Let
the  curve $\gamma$ in the affine plane $A^2$ be
defined in Cartesian coordinates by
$$x(t)=(\cos(2t)+5)\cos t, \hskip .7cm y(t)=(\cos(2t)+5)\sin t.$$
This curve is convex and has exactly two sextactic points and two
inflection points (which are
evidently not true inflection points), showing at the same time that
Theorem 1.2 (ii) and Theorem
1.3 are optimal. The affine curvature goes to negative infinity as one
approaches the inflection
points and has a local maximum between the inflection points. There is a
sketch of the curve in Figure B.2.
\medskip
\def\B{\vtop{\halign{\hfill##\hfill\cr
               {\epsfxsize=2in \epsfbox{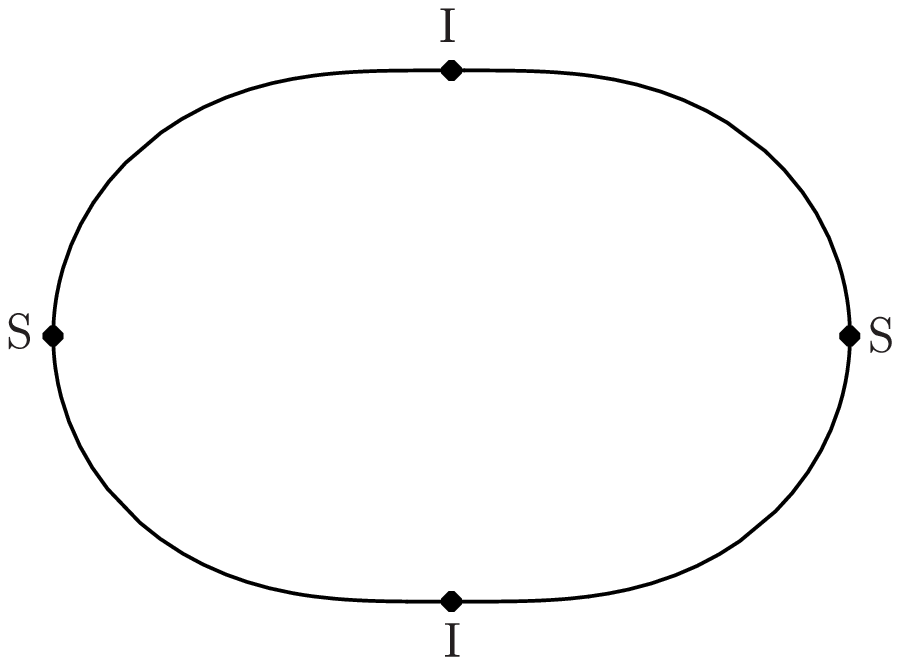}} \cr
            \noalign{\smallskip}{Figure B.2} \cr}}}
\centerline{\B}
\bigskip
\bigskip

\ni{\bf B.3 Example.} We consider the curve $\gamma$ in the affine plane
$A^2$ given in Cartesian
coordinates by
$$x(t)=(3+2\cos t)\cos t,\hskip .7cm y(t)=(3+2\cos t)\sin t.$$
This curve is not convex. It has two inflection points (both of them true)
and four sextactic
points. It shows that the last claim in Theorem 5.5 is optimal. There is a
sketch of the curve in Figure B.3. Notice that
we twice mark  $I$ and $S$ at the same place, since
the inflection points at $t\approx \pi\pm 0.352$ are
so close to the sextactic points at
$t \approx \pi\pm 0.335$ that one cannot distinguish between then in the
figure.
\medskip
\def\C{\vtop{\halign{\hfill##\hfill\cr
               {\epsfxsize=1.7in \epsfbox{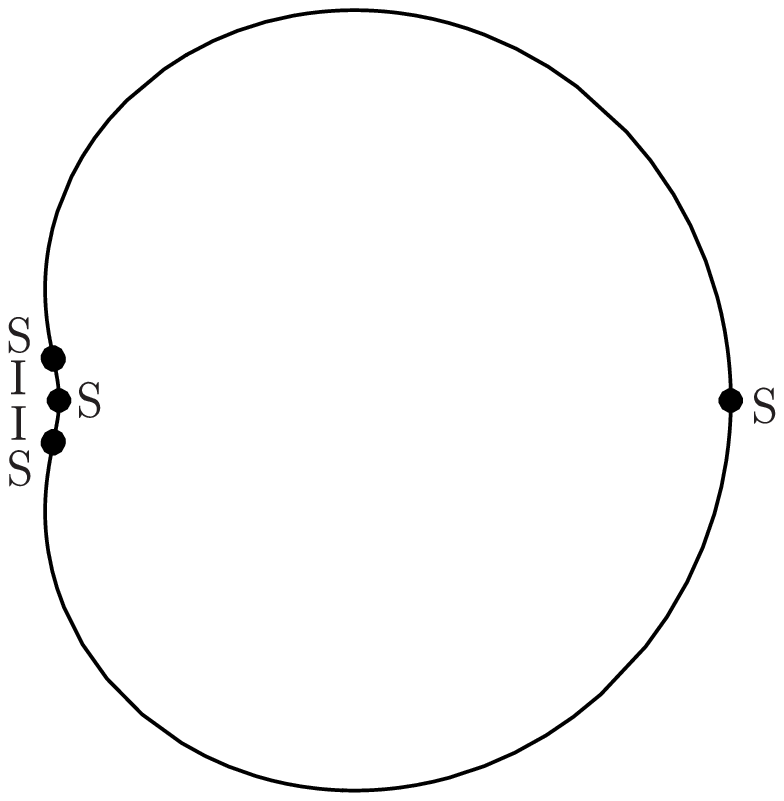}} \cr
            \noalign{\smallskip}{Figure B.3} \cr}}}
\centerline{\C}
\bigskip
\bigskip

\ni{\bf B.4 Example.}
This example  due to Izumiya and Sano [IS] also
shows that Theorem 1.3 is optimal. Here $\gamma$ is given in Cartesian
coordinates by
$$x(t)=(2+\cos t)\cos t, \hskip .7cm y(t)=(2+\cos t)\sin t.$$
The curve is convex. It has one inflection point and three sextactic
points. The affine curvature of
$\gamma$ has two local
maxima, one local minimum and it goes to negative infinity as one
approaches the inflection
point.  This example show that Theorem 1.3 and the last claim in Theorem
5.4 are optimal. There is a sketch of the curve in Figure B.4.
\medskip
\def\D{\vtop{\halign{\hfill##\hfill\cr
               {\epsfxsize=1.7in \epsfbox{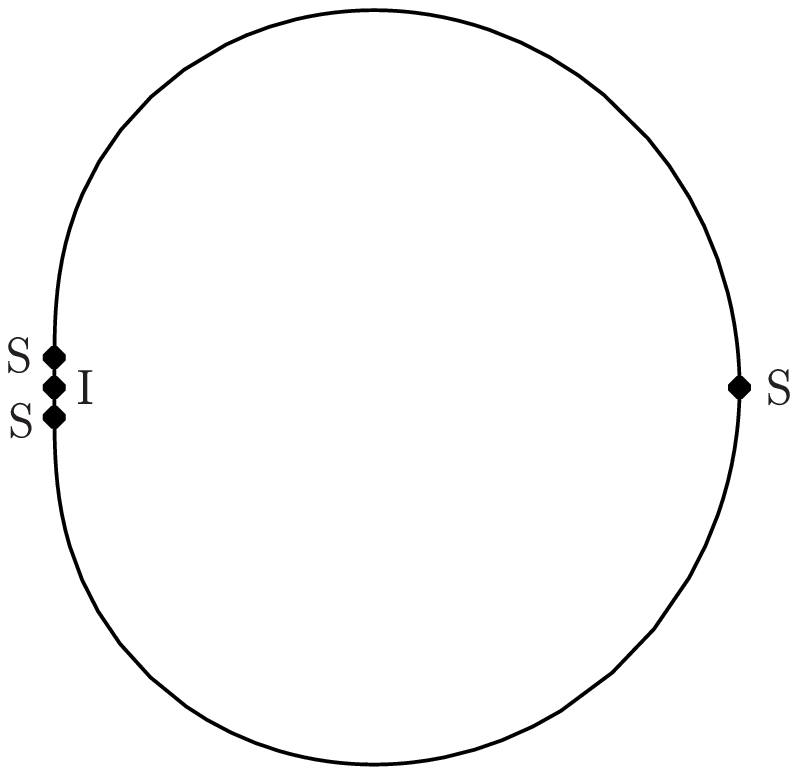}} \cr
            \noalign{\smallskip}{Figure B.4} \cr}}}
\centerline{\D}
\vfill
\eject

\ni{\bf  Appendix C: Sextactic points on a complex plane algebraic
curve}\medskip

We will give a proof of the theorem of Cayley [Ca 2] mentioned in the
introduction in
this  appendix using the theory of inflection points of linear systems as
explained in the
textbook [Mi], in which one can find explanations of all concepts used
here. The proof is
essentially only an adaptation of the methods used to prove the formula of
Pl\"ucker on the
number of inflection points in [Mi], p.~241. We also get a formula for the
number of
sextactic points when the inflection points are not all simple. This
probably also
follows from Cayley's method, but we believe that the methods below are
simpler. Other proofs of this theorem of Cayley can for example be found in
[Bt], [Bs], and
[Vi] together with references to further papers on the subject.

\medskip

We will be
considering the linear system of intersection divisors of conics. This
linear system  corresponds to
the Veronese embedding of the curve into $P^5({\bf C})$. We will therefore
really be studying the
number of points of higher order contact between a curve in $P^5({\bf C})$
and its osculating
hyperplanes. Such an approach was used by Barner in [Ba] to prove
Mukhopadhyaya's Theorem, see
also [Ar 3]. As can be seen in these papers, the method can also be used to
find the existence of
what is called an extatic point of a curve in $P^2$, i.e., the analogues of
sextactic points when
the conics are replaced by algebraic curves of some fixed degree, see [Ar
3], but very strong
conditions on the curve are needed. The number of extatic points of some
given order on an
algebraic curve in
$P^2({\bf C})$ can in principle also be determined as in the following
proof.

\medskip

\medskip
\ni{\bf C.1 Theorem.}  {\it Let $\gamma$ be a regular algebraic curve of
degree
$d$ in $P^2({\bf C})$. Then $\gamma$ has exactly $3d(4d-9)$ sextactic
points counted with
multiplicities if all inflection points of $\gamma$ are simple. If $\gamma$
has $k$ inflection
points with multiplicities $\nu_1,\dots,\nu_k$ respectively, then $\gamma$
has
$$3d(5d-11)-\sum_{i=1}^k 4\nu_i-3$$
sextactic points counted with multiplicities.}

\medskip
\ni{\it Proof.}  Here a conic will not be assumed to be regular. Let
$C$ be a conic. Then
$C$ induces a divisor
${\rm div}(C)$ on $\gamma$ by associating to $p\in\gamma$ the intersection
multiplicity of $C$
and $\gamma$ in $p$. By B\'ezout's theorem, the degree of ${\rm div}(C)$ is
equal to $2d$. The
collection of these divisors is a complete linear system
$Q$ of dimension $5$, i.e., $Q$ is a $g^5_{2d}$.

We have to determine the gap numbers for $Q$ at a point $p$. These are the
integers $\ell$ at
which the dimensions of the spaces in the sequence
$$
Q\supset Q(-p)\supset\dots\supset Q(-\ell p)\supset\dots
$$
change. We review that $Q(-\ell p)$ is the space of divisors in $Q$ that
meet $\gamma$
in $p$ with multiplicity $\ell$ at least.

Let us first assume that  $p$ is not an inflection point and that the
multiplicity with
which the osculating conic at $p$ meets $\gamma$ in $p$ is $\mu$. Then
$\dim Q(-\ell
p)=5-\ell$ for $\ell=1, 2, 3, 4, 5$, $\dim Q(-\ell p)=0$ for $\ell=5,\dots,
\mu$ and $
Q(-(\mu+1)p)=\emptyset$. Hence the gap sequence is $n_1=1$, $n_2=2$,
$n_3=3$, $n_4=4$, $n_5=5$
and $n_6=\mu+1$ if $\mu>5$.

The inflectionary weight of $p$ is by definition equal to
$$
w_p(Q)=\sum_{i=1}^6 (n_i-i).
$$
Hence $w_p(Q)=\mu-5$ if $p$ is not an inflection point. Notice that
$\mu-5$ is equal to
$0$ if $p$ is not an sextactic point. Otherwise $\mu-5$ is the multiplicity
of the sextactic
point.

Now let us assume that $p$ is an inflection point of $\gamma$ in which the
tangent line at
$\gamma$ and $\gamma$ meet with multiplicity $\mu$. The dimensions of
$Q(-p)$ and $Q(-2)$ do
not depend on whether we are at an inflection point or not, i.e., $\dim
Q(-p)=4$ and $\dim
Q(-2p)=3$. The spaces $Q(-3p)=\dots=Q(-\mu p)$ consist of the divisors of
conics that
are two lines, one of which is the tangent line, the other one arbitrary.
Hence $\dim Q(-\ell
p)=2$ for
$\ell=3,\dots,\mu$. The space $Q(-(\mu+1)p)$ consists of the divisors of
conics that are two
lines, one of which is the tangent line, the other passing through $p$.
Hence $\dim
Q(-(\mu+1)p)=1$. The spaces $Q(-(\mu+2)p=\dots=Q(-2\mu p)$ consist only of
the divisor of the
double tangent line. Hence $\dim Q(-(\mu+2)p)=\dots=\dim Q(-2\mu p)=0$. The
space $Q(-\ell
p)=\emptyset$ for $\ell \ge 2\mu+1$.
It follows that $n_1=1$, $n_2=2$, $n_3=3$, $n_4=\mu+1$, $n_5=\mu+2$,
$n_6=2\mu+1$. We
therefore have $w_p(Q)=(\mu-3)+(\mu-3)+(2\mu-5)=4\mu-11$. The multiplicity
of $p$ as an
inflection point is $\nu=\mu-2$. Hence $w_p(Q)=4\nu-3$. If $p$ is a simple
inflection point,
i.e., $\nu=1$, then $w_p(Q)=1$.

We are now going to use the formula
$$
\sum_{p\in\gamma}w_p(Q)=6(2d+5g-5),
$$
see [Mi], p.~241, where $g$ is the genus of $\gamma$, i.e.,
$g=(d-1)(d-2)/2$ by the Pl\"ucker
formula. Hence the number of sextactic points counted with multiplicities
is equal to
$$
3d(5d-11)-\sum_{i=1}^k 4\nu_i-3.
$$
If all inflection points are simple, i.e., $\nu_i=1$ for all $i$, then the
sum is equal to the
number of inflection points, which we know to be $3d(d-2)$. Hence the
number of sextactic
points is equal to $3d(4d-9)$ in that case, and we have finished the proof
of the theorem.
\qed

\medskip
\ni{\bf C.2 Example.} In this example we will explain the distribution of
inflection and sextactic
points on regular real and complex cubics. We have referred to the real
cubic in this
paper as an example for certain of our theorems being optimal.

We first consider the complex case. Let $\gamma$ be a regular complex plane
cubic.
First notice that a line and a cubic meet in three
points and a conic and cubic in six points. It follows that all inflection
and sextactic
points on $\gamma$ are simple. We therefore have precisely nine inflection
points and
precisely twenty seven sextactic points on $\gamma$. The distribution of
the inflection points is
well known.  If we choose one of the inflection points as the origin in the
group law of
the cubic $\gamma$ and denote it by $0$, then  a $p\in\gamma$ is an
inflection point if and only
if $3p=0$. Now one can show that all points $p\in \gamma$ with $2p=0$ are
sextactic. These are
not all sextactic points. In fact one can show that a point $p\in \gamma$
is either an
inflection or a sextactic point if and only if $6p=0$. Bearing in mind that
$\gamma$ as a
group is isomorphic to a torus ${\bf C}/\Lambda$, we see that the equation
$6p=0$ has thirty
six solutions as should be the case.

Now we come to the real parts of regular complex cubics.  A real cubic can
contain one or
two  branches. If it consists of one branch, it must correspond to the real
part of  ${\bf
C}/\Lambda$ and we see from the above description that it has precisely
three
inflection points and three sextactic
points. If the real cubic consists of two branches, one part must be the
real part of ${\bf
C}/\Lambda$, the other will be the image in ${\bf C}/\Lambda$ of the line
parallel to the real
axis passing through the center of the fundamental domain. Notice that this
second branch does not contain any inflection points and is therefore
strictly convex. Notice
also that it contains precisely six sextactic points.

\bigskip

\bigskip\bigskip
\ni{\bf Acknowledgment.} The authors wish to thank S. Izumiya, T. Sano
and T. Kurose for fruitful discussions.

\bigskip\bigskip\bigskip
\noindent {\bf References}\medskip

\item{[Ar 1]} V.I Arnold: {\it A ramified covering of ${\bf C}P^2\to S^4$,
hyperbolicity and
projective topology} (Russian). Russian original in: Sib. Mat. Zh. {\bf 29}
(1988), 36--47. English
translation in: Sib. Math. J.~{\bf 29} (1988), 717-726.

\item{[Ar 2]} V.I. Arnold: {\it Topological Invariants of Plane Curves and
Caustics.} University Lecture Series {\bf 5}. American Mathematical Society,
Providence, Rhode Island, 1994.

\item{[Ar 3]} V. Arnold: {\it
Remarks on the extatic points of plane curves.
}In: The Gelfand
Mathematical Seminars, 1993--1995, 11--22,
Birkh\"auser, Boston,
1996.

\item{[Ba]} M. Barner:
{\it \"Uber die Mindestanzahl station\"arer Schmiegebenen bei
geschlos\-senen strengkonvexen
Raumkurven.}
 Abh.~Math.~Sem.~Univ.~Hamburg {\bf 20} (1956), 196-215.

\item{[Bs]} A.B. Basset: {\it On sextactic and allied conics.} Quart.~J.
{\bf 46} (1915),
247--252.

\item{[Bt]} G. Battaglini: {Sui punti sestatici di una curva qualunque.}
Atti R.~Acc.~Lincei, Rend.~(Serie quarta)  {\bf IV}$_2$ (1888), 238-246.

\item{[Bl 1]} W. Blaschke: {\it \"Uber affine Geometrie VIII:
Die
Mindestzahl der sextaktischen Punkte einer Eilinie.} Leipziger
Berichte
{\bf 69} (1917), 321--324. Also in: Gesammelte Werke, Band
4,
153--156. Thales Verlag, Essen, 1985.

\item{[Bl 2]} W. Blaschke:
{\it Vorlesungen \"uber Differentialgeometrie
 II, Affine
Differentialgeometrie.} Springer-Verlag, Berlin, 1923.

\item{[Bo]} G.
Bol: {\it Projektive Differentialgeometrie, 1. Teil.}
Vandenhoeck \&
Ruprecht, G\"ottin\-gen, 1950.

\item{[Ca 1]} A. Cayley: {\it On the conic of five-pointic contact at any
point
of a plane curve.} Philosophical Transactions of the Royal Society of London
{\bf CXLIX} (1859), 371--400. Also in: The Collected Mathematical Papers,
Vol.~IV,
Cambridge: at the University Press, 1891.

\item{[Ca 2]} A. Cayley: {\it On the sextactic points of a plane curve.}
Philosophical Transactions
of the Royal Society of London {\bf CLV} (1865), 548--578. Also in: The
Collected Mathematical
Papers, Vol.~V, Cambridge: at the University Press, 1892.

\item{[Fa]} Fr.~Fabricius-Bjerre: {\it On a conjecture of G. Bol.}
Math.~Scand.~{\bf 40}
(1977), 194--196.

\item{[GMO]}  L. Guieu, E. Mourre \&
V. Yu. Ovsienko: {\it  Theorem on
six vertices of a plane curve via
Sturm theory.} In: The Arnold-Gelfand
Mathematical Seminars, 257--266,
Birkh\"auser, Boston, 1997.

\item{[IS]} S. Izumiya \& T. Sano: {\it Private Communication}. 1998.

\item{[Kn]} H. Kneser: {\it Neuer Beweis
des Vierscheitelsatzes.}
Christiaan Huygens {\bf 2} (1922/23),
315--318.

\item{[Mi]} R. Miranda: {\it Algebraic curves and Riemann surfaces.}
Graduate Studies in
Mathematics {\bf 5}. American Mathematical Society,
Providence, Rhode Island, 1995.

\item{[M\"o]} A. F. M\"obius: {\it \"Uber die Grundformen
der Linien
der dritten Ordnung.} Abhandlungen
der
K\"onigl.~S\"achs.~Gesellschaft der Wissenschaften, math.-phys.
Klasse I (1852), 1--82. Also in:  Gesammelte Werke, vol.
II, Verlag
von S. Hirzel, Leipzig, 1886, 89--176.

\item{[Mu 1]}
S. Mukhopadhyaya: {\it New methods in the geometry of a plane arc, I.}
Bull. Calcutta Math. Soc. {\bf 1} (1909), 31-- 37. Also in:
Collected geometrical papers, vol. I.   Calcutta University Press,
Calcutta, 1929, 13--20.

\item{[Mu 2]} S. Mukhopadhyaya: {\it Sur les
nouvelles m\'ethodes de
g\'eometrie.} C. R. S\'eance Soc.~Math.
France, ann\'ee 1933 (1934),
41-45.

\item{[TU 1]} G. Thorbergsson \& M. Umehara: {\it A
unified approach to the
four vertex
theorems, II.} In: Differential  and symplectic  topology of
knots and curves (edited by S. Tabachnikov), Amer.~Math.~Soc.~Transl.~{\bf
190}
(1999), 229--252.

\item{[TU 2]} G. Thorbergsson \& M. Umehara: {\it On global properties of
flexes of periodic
functions.} Preprint 2000.

\item{[Um]} M.
Umehara: {\it A unified approach to the four vertex
theorems, I.}
In: Differential  and symplectic  topology of
knots and curves (edited by S. Tabachnikov), Amer.~Math.~Soc.~Transl.~{\bf
190}
(1999), 185--228.

\item{[Vi]} A.O. Viro: {\it Differential geometry \lq\lq in the large\rq\rq\
of plane
algebraic curves, and integral formulas for invariants of singularities}
(Russian).
Zap.~Nauchn. Sem.~S.-Petersburg. Otdel.~Mat.~Inst.~Steklov. (POMI) {\bf 231}
(1995), 255-268.

\item{[Wa]} C.T.C. Wall: {\it Duality of real projective plane curves:
Klein's equation.}
Topology {\bf 35} (1996), 355--362.
\vskip 1.in

\hbox{\parindent=0pt\parskip=0pt
\vbox{\hsize=2.7truein
\obeylines
{
Gudlaugur Thorbergsson
Mathematisches Institut
Universit\"at zu K\"oln
Weyertal 86-90
50931 K\"oln
Germany
}\medskip
gthorbergsson@mi.uni-koeln.de
}\hskip 1.5truecm
\vbox{\hsize=3.7truein
\obeylines
{
Masaaki Umehara
Department of Mathematics
Graduate School of Science,
Hiroshima University
Higashi-Hiroshima 739-8526
Japan}\medskip
umehara@math.sci-hiroshima-u.ac.jp
}
}
\bye